\documentclass[12pt]{amsart}
\usepackage{bbm}
\usepackage{amsfonts}
\usepackage{amssymb}
\usepackage{mathrsfs}
\usepackage{}
\usepackage{stmaryrd}
\usepackage{txfonts}
\usepackage{amssymb, amsmath}
\usepackage{graphicx}
\usepackage{graphics}
\usepackage{latexsym}
\usepackage{amsmath}
\usepackage{amssymb,amsthm,amsfonts}
\usepackage{amscd}
\usepackage[arrow, matrix, curve]{xy}
\usepackage{syntonly}
\usepackage{extarrows}
\theoremstyle{plain}
\newtheorem{thm}{Theorem}[section]
\newtheorem{theorem}[thm]{Theorem}
\newtheorem{lemma}[thm]{Lemma}

\theoremstyle{definition}

\newtheorem{remark}[thm]{Remark}

\newtheorem{definition}[thm]{Definition}

\newtheorem{example}[thm]{Example}

\begin{document}
\title{Kauffman polynomial from a generalized Yang-Yang function}
\author {Sen Hu and Peng Liu}
\maketitle

\begin{abstract}
For the fundamental representations of the simple Lie algebras of type $B_{n}$, $C_{n}$ and $D_{n}$, we derive the braiding and fusion matrices from the generalized Yang-Yang function and prove that the corresponding knot invariants  are Kauffman polynomial.
\end{abstract}

\tableofcontents
\section{Introduction}
After E.Witten's work \cite{W1} on Jones polynomial from quantum field theory, knots theory has inspired lots of interest of physicists and mathematicians. Recently D.Gaiotto and E.Witten have developed a new method to derive knot invariants from gauge theory. The free field realization of Virasoro conformal blocks was used to convert the integral of Chern-Simons functional over the infinite dimensional moduli space of connections into the integral of a Yang-Yang function over a finite dimensional parameter space. What's more important, Lefschetz thimbles naturally give a representation of the braiding operation and therefore provide a powerful tool to study the wall-crossing phenomena. In \cite{HL}, we used the generalized Yang-Yang function to study the braiding matrix and knot invariants for $A_{n}$ Lie algebras. In this paper, we use quantum mechanics to describe this method, derive the braiding and fusion matrices for the simple Lie algebras of type $B_{n},C_{n}$ and $D_{n}$ and prove that the corresponding knot invariants are Kauffman polynomial.

In section 2, we briefly review getting braiding formula from Lefschetz thimbles of the generalized Yang-Yang functions for both cases of
symmetry breaking and without symmetry breaking. This part of results was developed in the previous paper \cite{HL}. 

In section 3, we define the braiding and fusion operators in quantum mechanics language. Then,  for the simple Lie algebras of type $B_{n},C_{n}$ and $D_{n}$, we derive the braiding and fusion matrices  and several relations between them from thimbles of generalized Yang-Yang functions in section 3. 

In section 4, we prove that the corresponding knot invariants are Kauffman polynomial.

\section{Braiding formula from Lefschetz thimbles of the generalized Yang-Yang function}

The real part of a holomorphic function on a Hermitian manifold, as a Morse function, has some nice properties.
\begin{lemma}\label{lem:thimble}
For a holomorphic function on a Hermitian manifold $M$ $\mathrm{dim}_{\mathbb{R}}M=2d$, if its real part $h$ is a Morse function on $M$ (i.e. the Hessian matrix of $h$ is non-degenerate and its critical points are isolated), then $(1)$ the gradient flow of $h$ keeps the imaginary part invariant; $(2)$ the index of each critical point of $h$ is $d$.
\end{lemma}

From this lemma, Lefschetz thimbles are naturally defined:
\begin{definition}
Assuming that $h$ defined as in the last lemma is a Morse function on $M$, a cycle is called a Lefschetz thimble associated to $I$, denoted by $\mathcal {J}$, if all points of it can be reached by the gradient flow of $h$ starting from a critical point $I$.
\end{definition}
Lefschetz thimbles are middle dimensional cycles of $M$.
Here is an example:
\begin{example}
$f(x)=i\lambda ( \dfrac {x^{3}} {3}-x), x\in \mathbb{C}$ is called Airy function \cite{W2}, where $\lambda=a+bi$ is a complex constant with $b>0$.
The critical points of the holomorphic function $f(x)$ are $x=\pm 1$, denoted by $P_{\pm}=\pm 1$. $\mathrm{Im} f(P_{+})=-\frac{2}{3}a$ and $\mathrm{Im} f(P_{-})=\frac{2}{3}a$. $\mathrm{Im} f(P_{+})=\mathrm{Im} f(P_{-})$ if and only if $a=0$. Thus, from Lemma \ref{lem:thimble}, there is a gradient flow defined by Morse function $Re f(x)$ connecting $P_{+}$ with $P_{-}$ if and only if $a=0$. When $a=0$, $f(x)=-b(\dfrac {x^{3}} {3}-x)$. The gradient flow connecting $P_{+}$ with $P_{-}$ is on the real axis of $x$ plane and the imaginary part of $f$ is zero along this flow. When $a\neq0$, there is no gradient flow connecting $P_{+}$ with $P_{-}$. As is shown in Figure \ref{fig:thimble}, the picture (a), (b) and (c) describe the gradient flows starting from $P_{+}$ and $P_{-}$ with $a=1$, $a=0$ and $a=-1$ respectively.

\begin{figure}
\begin{center}
\includegraphics[width=12cm,height=6cm]{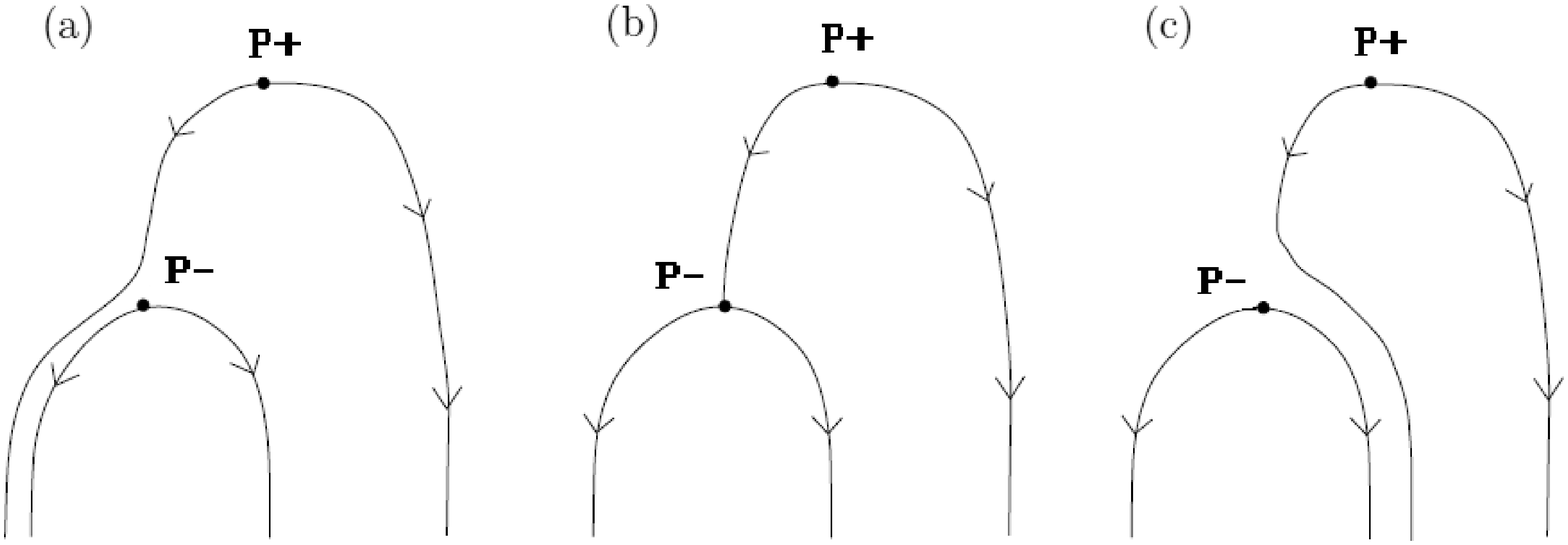}
\caption{\label{fig:thimble} Wall-crossing}
\end{center}
\end{figure}

 If $a$ is continuously changed from $1$ to $-1$ on the $\lambda$ plane with $b>0$, the thimble $\mathcal{J}_{+}$ associated to $P_{+}$ and the thimble $\mathcal{J}_{-}$ associated to $P_{-}$ will be transformed into $\mathcal{J}_{+}^{'}$ and  $\mathcal{J}_{-}^{'}$:\begin{equation}\left(
                                                                               \begin{array}{c}
                                                                                 \mathcal{J}_{+}^{'} \\
                                                                                 \mathcal{J}_{-}^{'} \\
                                                                               \end{array}
                                                                             \right)=\left(
                                                                                       \begin{array}{cc}
                                                                                         1 & \pm 1 \\
                                                                                         0 & 1 \\
                                                                                       \end{array}
                                                                                     \right)\left(
                                                                               \begin{array}{c}
                                                                                 \mathcal{J}_{+} \\
                                                                                 \mathcal{J}_{-} \\
                                                                               \end{array}
                                                                             \right).
                                                                             \end{equation}
$\mathcal{J}_{-}$ is unchanged, but $\mathcal{J}_{+}$ will receive an additional term. This phenomena also appears when we change $a$ from $1$ to $-1$ continuously with $b<0$ .
When $a=0$, two rays $b>0$ and $b<0$ on the $\lambda$ plane are called Stokes rays (Stokes walls). Passing through a Stokes ray is called wall-crossing. When wall-crossing happens, the thimble will receive an additional term.
\end{example}
Considering distinct points $z_{1},...,z_{d}$ on $\mathbb{C}$, we associate each point with an irreducible highest weight representation $V_{\lambda_{a}}$ of $\mathfrak{g}$, where $\lambda_{a}$ is a dominant integral weight. Thus $V_{\lambda_{a}}$ is a finite dimensional irreducible highest weight representation of $\mathfrak{g}$.  $V_{(\lambda_{a})}\triangleq V_{\lambda_{1}}\otimes V_{\lambda_{2}}\otimes...\otimes V_{\lambda_{d}} $. $\Pi =\{\alpha_{1},\alpha_{2},...,\alpha_{n}\}$ is the set of the simple roots of $\mathfrak{g}$.
$w_{j} (j=1,2,...,\mathbbm{q})$ are distinct points on $\mathbb{C}$ different from $z_{a}$. Each $w_{j}$ is associated with a simple root $\alpha_{i_{j}}$ of $\mathfrak{g}$, where $i_{j}\in \{1,2,...,n\}$.

\begin{equation}
\mathscr{W}(w,z)=\sum _{j,a}(\alpha_{i_{j}},\lambda _{a})\ln( w_{j}-z_{a}) -\sum _{j< s}(\alpha_{i_{j}},\alpha_{i_{s}})\ln ( w_{j}-w_{s})-\sum _{a< b}(\lambda_{a},\lambda_{b})\ln(z_{a}-z_{b})
\end{equation}
is the generalized Yang-Yang function (see \cite{FFR}, \cite{ATY} and \cite{F} for more details) associated to the representation $V_{\lambda_{a}}$ of $\mathfrak{g}$.

We consider the integration $\int _{\mathcal{J}_{\mathbbm{q}}}e^{-\frac{\mathscr{W}}{k+h^{\vee }}}\prod _{j}dw_{j},$ where $\mathcal{J}_{\mathbbm{q}}$ is the thimble of the generalized Yang-Yang function with $\mathbbm{q}$ variables $w_{1}, w_{2}, ..., w_{\mathbbm{q}}$. It is a sub-manifold with $\mathrm{dim}_{\mathbb{R}}\mathcal{J}_{\mathbbm{q}}=\mathbbm{q}$ in $\mathbb{C}^{\mathbbm{q}}$.  Here we only consider the case $d=2$ and $\lambda_{1}=\lambda_{2}$. Without losing generality, we assume that $z_{1}$ and $z_{2}$ have the same real part and $\mathrm{Im} z_{1}>\mathrm{Im} z_{2}$. Then we rotate $z_{1}$ and $z_{2}$ clockwise by $\pi$ around the middle point of them. The multiple valued Yang-Yang function will produce an phase factor under this braiding transformation.
It can be easily obtained by the following lemma.
\begin{lemma}\label{lem:img part inv}
$f(w_{j},z_{1},z_{2})$ is a holomorphic function of $w_{j}$ with two complex parameters $z_{1}$ and $z_{2}$. If $\Gamma$ is a thimble associated to the real part of $f(w_{j},z_{1},z_{2})$, then the phase factor of the integral $\int _{\Gamma }e^{f(w_{j},z_{1},z_{2})}\prod _{j}dw_{j}$ coming from the braiding without wall-crossing is equal to the phase factor of $e^{f(w_{c},z_{1},z_{2})}$ under the braiding, where $w_c$ is the the critical point of $f(w_{j},z_{1},z_{2})$.
\end{lemma}

From the Lemma \ref{lem:img part inv}, we have
 \begin{equation}\label{braiding formula in NSB}
\mathcal{B}\mathcal{J}_{\mathbbm{q}}=(-1)^{\mathbbm{q}}q^{-\frac{1}{2}[(\lambda_{1},\lambda_{2})+\sum_{j<s}(\alpha_{i_{j}},\alpha_{i_{s}})-\sum_{j,a}(\alpha_{i_{j}},\lambda_{a})]}\mathcal{J}_{\mathbbm{q}}.
   \end{equation}
The $(-1)^{\mathbbm{q}}$ comes from the fact that the braiding changes the direction of each dimension into the opposite direction and the thimble $\mathcal{J}_{\mathbbm{q}}$ is $\mathbbm{q}$ dimensional. We formally define $\mathcal{J}_{0}$ to record the phase factor of $e^{-\frac{\mathscr{W}(z_{1},z_{2})}{k+h^{\vee }}}$:
\begin{equation}\label{J0}
\mathcal{B}\mathcal{J}_{0}=q^{-\frac{1}{2}(\lambda_{1},\lambda_{2})}\mathcal{J}_{0}.
\end{equation}

With a positive real parameter $c$, \begin{multline}
\mathscr{W}(w_{j},z_{1},z_{2},c)
=\sum _{j}(\alpha_{i_{j}},\lambda _{1})\ln( w_{j}-z_{1})+ \sum _{j}(\alpha_{i_{j}},\lambda _{2})\ln( w_{j}-z_{2})-\sum _{j<s}(\alpha_{i_{j}},\alpha_{i_{s}})\ln ( w_{j}-w_{s})\\-(\lambda_{1},\lambda_{2})\ln(z_{1}-z_{2})-c\left(\sum_{j}w_{j}-\frac{1}{2}\sum_{a}\parallel \lambda_{a}\parallel z_{a}\right)
\end{multline} is called a generalized Yang-Yang function with symmetry breaking.
When $c\rightarrow 0$, it goes back to the generalized Yang-Yang function, i.e. $\underset{c\rightarrow 0}{\lim}\mathscr{W}(w_{j},z_{1},z_{2},c)=\mathscr{W}(w_{j},z_{1},z_{2}).$
When $c\rightarrow \infty$,
the critical point equation or Bethe equation
\begin{equation}\frac{\partial\mathscr{W}(w_{j},z_{1},z_{2},c)}{\partial w_{j}}=0,\quad j=1,2,...,\mathbbm{q}\end{equation}
has solutions:
$$w_{j}=z_{1}+o (\frac{1}{c})\hbox{ or }z_{2}+o(\frac{1}{c}).$$
Let \begin{equation}
S_{1}=\{j|w_{j}=z_{1}+o (\frac{1}{c})\}
\end{equation}
and
\begin{equation}
S_{2}=\{j|w_{j}=z_{2}+o (\frac{1}{c})\},
\end{equation}
then $$\#S_{1}+\#S_{2}=\mathbbm{q}.$$
Consider special solutions satisfying the following two conditions:\begin{enumerate}
                                                              \item $\#S_{1}\leq m$, $\#S_{2}\leq m$, where $m=\mathrm{dim} V_{\lambda}-1$;
                                                              \item $\lambda-\sum_{j\in S_{1}}\alpha_{i_{j}}$ and $\lambda-\sum_{j\in S_{2}}\alpha_{i_{j}}$ are weights of the representation $V_{\lambda}$.
                                                            \end{enumerate}
Denote the thimble associated to the critical point of this type as
$\mathcal{J}_{s,\mathbbm{q}-s}$, where $s=\#S_{1}$. We can formally define $\mathcal{J}_{0,0}=\mathcal{J}_{0}$ as in (\ref{J0}), then $0\leq\mathbbm{q}\leq2m$.
Now combining all thimbles together with respect to $\mathbbm{q}$ from $0$ to $2m$, we have totally
$(m+1)^{2}$ different thimbles.
These $(m+1)^{2}$ thimbles
associated to special solutions of the Bethe equation in the symmetry breaking case naturally form a set of basis
of the representation space $V_{\lambda_{1}}\otimes V_{\lambda_{2}}$. Each thimble $\mathcal{J}_{s,\mathbbm{q}-s}$ corresponds to a weight vector
with weight $\lambda_{1}-\sum_{j\in S_{1}}\alpha_{i_{j}},\lambda_{2}-\sum_{j\in S_{2}}\alpha_{i_{j}}$
in the representation space $V_{\lambda_{1}}\otimes V_{\lambda_{2}}$.
$$V_{\lambda_{1}}\otimes V_{\lambda_{2}}=\oplus_{\mathbbm{q}=0}^{2m}V_{\mathbbm{q}},
$$ where $V_{\mathbbm{q}}$ is a linear space generated by $\{v_{s}\otimes v_{\mathbbm{q}-s}, 0\leq s\leq\mathbbm{q} \mid v_{s}\otimes v_{\mathbbm{q}-s}\in V_{\lambda_{1}}\otimes V_{\lambda_{2}}, v_{s}\otimes v_{\mathbbm{q}-s} \hbox { is a vector with weight } (\lambda-\sum_{j\in S_{1}}\alpha_{i_{j}},\lambda-\sum_{j\in S_{2}}\alpha_{i_{j}})  \}\subseteq V_{\lambda_{1}}\otimes V_{\lambda_{2}}$. Clearly, the braiding does not change the dimension of the thimble. Therefore, $V_{\mathbbm{q}}$ is an invariant subspace of the braiding transformation.

The braiding of the thimble of the generalized Yang-Yang function with symmetry breaking is
\begin{equation}\label{braiding formula in SB}
\mathcal {B}\mathcal{J}_{s,\mathbbm{q}-s}=q^{-\frac{1}{2}(\lambda_{1}-\sum_{1\leqq j\leqq s}\alpha_{i_{j}},\lambda_{2}-\sum_{s+1\leqq j\leqq \mathbbm{q}}\alpha_{i_{j}})}\mathcal{J}_{\mathbbm{q}-s ,s}+ w.c.t.,
\end{equation}
where w.c.t represents unknown wall-crossing terms. See \cite{HL} for more details.

Two properties of wall-crossing should be noticed. First, wall-crossing in the braiding transformation does not create or annihilate any simple root, but only transfers
 them from one location to another.
We call this property the conservation law of wall-crossing.
If the total types and numbers of the simple roots of two thimbles are different,
the Yang-Yang functions of them are two different functions. From the definition of the braiding of the
thimble, the braiding transformation only appears between the thimbles from one holomorphic function. Thus this
property is natural. Second, the transfer of simple roots in the wall-crossing can
only be from $z_{2}$ to $z_{1}$. The gradient flows in the symmetry breaking case are from $z_{1}$ and
$z_{2}$ to the infinity in the positive direction of the real axis in $w$ plane. In our assumption, $z_{1}$ and $z_{2}$ have
 the same real part and $\mathrm{Im} z_{1}>\mathrm{Im} z_{2}$. Therefore, in the clockwise braiding, the wall-crossing appears when there is a gradient flow started from $z_{2}$ passing through $z_{1}$.
Thus the only possible transfer of simple roots is from $z_{2}$
 to $z_{1}$. These two properties tell us that if we choose proper basis, the braiding matrix will be a diagonal
partitioned matrix and each block in the diagonal will be a triangular matrix.
 Thus, they are actually sub-representations of the braiding.

\section{Quantum mechanics and knot invariants}
\begin{figure}
\begin{center}
\includegraphics[width=8cm,height=4cm]{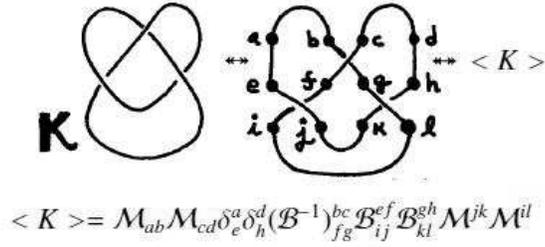}
\caption{\label{fig:knot}Decomposition of a link diagram}
\end{center}
\end{figure}

\begin{figure}
\begin{center}
\includegraphics[width=6cm,height=3cm]{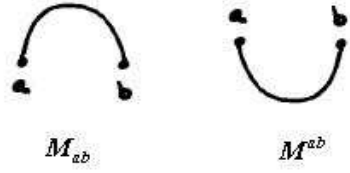}
\caption{\label{fig:AC} Annihilation $\mathcal{M}_{ab}$ and creation $\mathcal{M}^{ab}$}
\end{center}
\end{figure}

In \cite{K}, quantum mechanics was used to study knot invariants. We give a brief review of it in this section.
After the projection on a plane, knots can be decomposed as two strands braiding $\mathcal{B}$, its inverse $\mathcal{B}^{-1}$, annihilation $\mathcal{M}_{ab}$, creation $\mathcal{M}^{ab}$ and identities, as is shown in Figure \ref{fig:knot} and Figure \ref{fig:AC}
.

In quantum mechanics, the probability amplitude for the concatenation of processes is obtained by summing the products of the amplitudes  of the intermediate configurations in the process over all possible internal configurations. The initial states $v_{I}$, intermediate states $v_{M}$ and final states $v_{F}$ form vector spaces $V_{I}$, $V_{M}$ and $V_{F}$ respectively. Denote the first process by $P_{ab}$ from the initial state $v_{a}(a\in I)$ to the intermediate state $v_{b}(b\in M)$, then the $P_{ab}$ is a transformation from the $V_{I}$ to the $V_{M}$. The second process $Q_{cd}$ is from the intermediate state $v_{c}(c\in M)$ to the final state $v_{d}(d\in F)$. Then the amplitude of the concatenation processes starts from $v_{a}(a\in I)$ and ends up with $v_{d}(d\in F)$ is $\sum_{c}P_{ac}Q_{cd}=(PQ)_{ad}$. Knots can be thought as a process starting from and ending up with a vacuum state with braiding, its inverse, creation, annihilation and identities as intermediate configurations. Thus knot invariants are vacuum expectations of a quantum mechanics system. For example, in Figure \ref{fig:knot}, the invariant $<K>$ of the knot $K$  is $$<K>=\mathcal{M}_{ab}\mathcal{M}_{cd}\delta_{e}^{a}\delta_{h}^{d}(\mathcal{B}^{-1})_{fg}^{bc}\mathcal{B}_{ij}^{ef}\mathcal{B}_{kl}^{gh}\mathcal{M}^{jk}\mathcal{M}^{il},$$
where we use Einstein notation for summation.

Now we focus our attention on the state space and the operators acting on it. To derive  knot invariants from a general simple Lie algebra, we define the state space as the representation space $V_{\lambda}$ of the finite dimensional irreducible highest weight representation of the simple Lie algebra $\mathfrak{g}$ with the highest weight $\lambda$. Then the braiding operator is defined as $\mathcal {B}:V_{\lambda_{1}}\otimes V_{\lambda_{2}}\longrightarrow V_{\lambda_{1}}\otimes V_{\lambda_{2}}$. We assume $\mathrm{dim} V_{\lambda}=m+1$. Therefore, the braiding matrix is an $(m+1)^{2}\times (m+1)^{2}$ matrix.

Annihilator is a function $f:V_{\lambda_{1}}\otimes V_{\lambda_{2}}\longrightarrow \mathbb{C}$. It is equivalent to $\mathcal{M}:V_{\lambda_{1}}\longrightarrow V_{\lambda_{2}}^{*}\cong V_{\lambda_{1}}$
and can be represented as an $(m+1)\times (m+1)$ matrix. Creator is a map $g:\mathbb{C}\longrightarrow V_{\lambda_{1}}\otimes V_{\lambda_{2}}$.
It can also be represented as an $(m+1)\times (m+1)$ matrix. We denote the amplitudes for annihilation and creation between two states $v_{a}$ and $v_{b}$ as $\mathcal{M}_{ab}$ and $\mathcal{M}^{ab}$ respectively in Figure \ref{fig:AC}.
\begin{figure}
\begin{center}
\includegraphics[width=8cm,height=4cm]{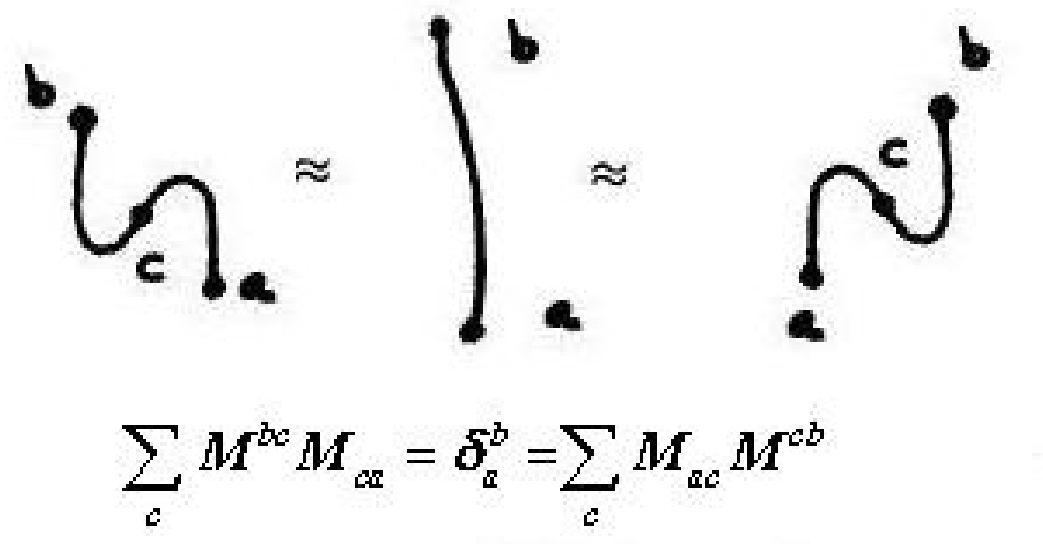}
\caption{\label{fig:TM}Invariance under the topological moves}
\end{center}
\end{figure}
As is shown in Figure \ref{fig:TM}, to be invariant under the topological moves, they should be inverse to each other:
\begin{equation}
\sum_{c}\mathcal{M}^{bc}\mathcal{M}_{ca}=\sum_{c}\mathcal{M}_{ac}\mathcal{M}^{cb}=\delta ^{b}_{a}.
\end{equation}
We call the matrix of amplitudes for the annihilation of two strands fusion matrix $\mathcal{M}$, then the matrix of amplitudes for the creation of two strands is just its inverse $\mathcal{M}^{-1}$.
For representations $V_{\lambda_{1}}$ and $V_{\lambda_{2}}$ ($\lambda_{1}=\lambda_{2}$) mentioned above, we define two vectors $v_{s}$ of weight $\lambda-\alpha_{i_{1}}-...-\alpha_{i_{s}}$ and $v_{m-s}$ of weight $\lambda-\alpha_{j_{1}}-...-\alpha_{j_{m-s}}$ in $V_{\lambda_{1}}$ and $V_{\lambda_{2}}$ to be complementary to each other. Since only two complementary states can fuse into a vacuum state,
fusion amplitudes are nonzero only between two complementary vectors.  Then fusion matrix can be written as
$\mathcal{M}:V_{\lambda}\longrightarrow V_{\lambda}$
            \begin{equation}\mathcal{M}\left(
                                                                            \begin{array}{c}
                                                                              v_{m}\\
                                                                              v_{m-1}\\
                                                                              ... \\
                                                                              v_{1} \\
                                                                              v_{0} \\
                                                                            \end{array}
                                                                          \right)
            =\left(
                                                 \begin{array}{ccccc}
                                                    &  &  &  &  \mathcal{M}_{m0}\\
                                                    & 0 &  & \mathcal{M}_{m-11} &  \\
                                                    &  & ... &  &  \\
                                                    & \mathcal{M}_{1m-1} &  & 0 &  \\
                                                   \mathcal{M}_{0m} &  &  &  &  \\
                                                 \end{array}
                                               \right)\left(
                                                                            \begin{array}{c}
                                                                              v_{m}\\
                                                                              v_{m-1}\\
                                                                              ... \\
                                                                              v_{1} \\
                                                                              v_{0} \\
                                                                            \end{array}
                                                                          \right).
            \end{equation}

\section{Braiding and fusion matrices for simple Lie algebras $B_{n},C_{n}$ and $D_{n}$}

\subsection{Fundamental representation of $B_{n}$}

\begin{figure}
\begin{center}
\includegraphics[width=2cm,height=10cm]{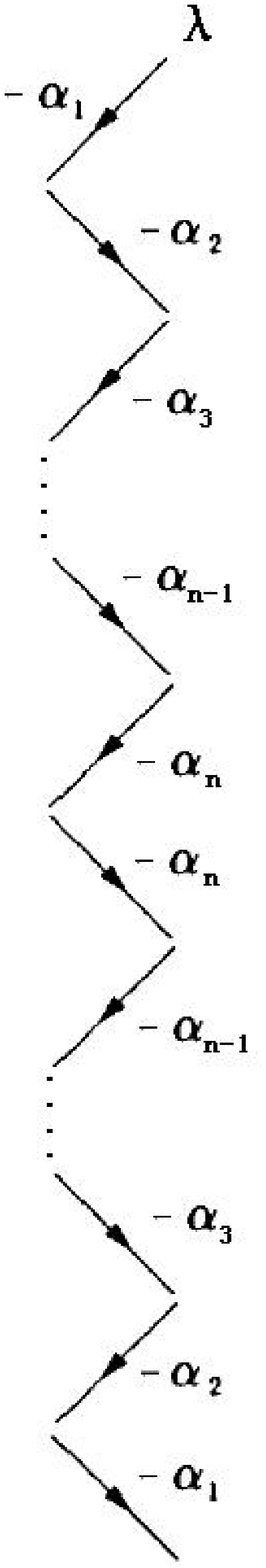}
\caption{\label{fig:Bn} Weights for the fundamental representation of $B_{n}$}
\end{center}
\end{figure}

$B_{n}$ has $n$ simple roots $\alpha_{i}, i=1,2,...,n$. Cartan matrix of $B_{n}$ is$$\left(
                                                     \begin{array}{ccccc}
                                                       2 & -1 &   &   &   \\
                                                       -1 & 2 & -1 &   &   \\
                                                         & -1 & ... &   &   \\
                                                         &   &   & 2 & -2 \\
                                                         &   &   & -1 & 2 \\
                                                     \end{array}
                                                   \right).$$
The highest weight of the fundamental representation $V_{\lambda}$ of $B_{n}$ is $\lambda=(1,0,...,0)$. Here we use Dynkin label: weights are represented in the fundamental weight basis. They are  $\lambda,\lambda-\alpha_{1},\lambda-\alpha_{1}-\alpha_{2},...,\lambda-\sum_{i=1}^{n}\alpha_{i}, \lambda-\sum_{i=1}^{n}\alpha_{i}-\alpha_{n},\lambda-\sum_{i=1}^{n}\alpha_{i}-\alpha_{n}-\alpha_{n-1},...,\lambda-\sum_{i=1}^{n}2\alpha_{i}$ , as is shown in Figure  \ref{fig:Bn}. It is a $2n+1$ dimensional representation and naturally gives an order on weights of the fundamental representation. The ordered weights are denoted by $\lambda^{0},\lambda^{1},...,\lambda^{2n}$.

\begin{lemma}\label{lem:Bn inner prd}
\begin{equation}
(\lambda^{s},\lambda^{t})=\left\{
                            \begin{array}{ll}
                              1, & \hbox{$s+t\neq2n,s=t$;} \\
                              0, & \hbox{$s+t\neq2n,s\neq t$;} \\
                              0, & \hbox{$s+t=2n,s=t $;} \\
                              -1, & \hbox{$s+t=2n,s\neq t $,}
                            \end{array}
                          \right.
\end{equation}where $s,t=0,1,...,2n$.
\end{lemma}
Proof: The proof is straightforward.

It should be noticed that these inner products are irrelative to the rank $n$.

The fundamental representation of $B_{n}$ has the following duality property.
\begin{lemma}
\begin{equation}
(\lambda^{s},\lambda^{t})=(\lambda^{2n-s},\lambda^{2n-t}),
\end{equation} where $s,t=0,1,...2n$.
\end{lemma}
Proof: $2n-s+2n-t=2n,2n-s=2n-t$ if and only if $s+t=2n,s=t$. $2n-s+2n-t\neq2n,2n-s=2n-t$ if and only if $s+t\neq2n,s=t$. From Lemma \ref{lem:Bn inner prd}, the proof is straightforward.

We use the equation (\ref{braiding formula in NSB}) to compute the phase factor coming from the braiding of the thimble $\mathcal{J}_{2n}$ of the generalized Yang-Yang function without symmetry breaking.
In the fundamental representation of $B_{n}$ Lie algebra, the simple roots $\alpha_{i_{j}}$ are $\alpha_{i_{j}}=\left\{
                                                                            \begin{array}{ll}
                                                                              \alpha_{j}, & \hbox{$1\leq j\leq n$;} \\
                                                                              \alpha_{2n+1-j}, & \hbox{$n<j\leq 2n$.}
                                                                            \end{array}
                                                                          \right.$

\begin{figure}
\begin{center}
\includegraphics[width=8cm,height=4cm]{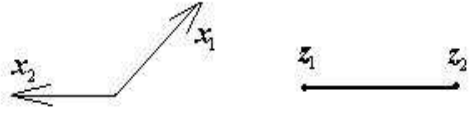}
\caption{\label{fig:BB}$\mathcal{J}_{2n}$ before braiding}
\end{center}
\end{figure}

\begin{figure}
\begin{center}
\includegraphics[width=8cm,height=4cm]{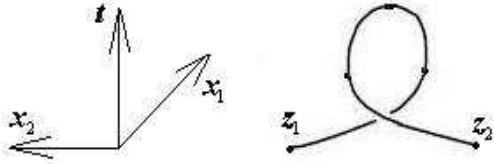}
\caption{\label{fig:AB}$\mathcal{J}_{2n}$ after braiding and projection}
\end{center}
\end{figure}

Thus \begin{equation}\label{braiding of B in NSB}
\mathcal{B}\mathcal{J}_{2n}=(-1)^{2n}q^{-\frac{1}{2}[(\lambda_{1},\lambda_{2})+\sum_{j<s}(\alpha_{i_{j}},\alpha_{i_{s}})-\sum_{j,\alpha}(\alpha_{i_{j}},\lambda_{\alpha})]}\mathcal{J}_{2n}=q^{n}\mathcal{J}_{2n}.
\end{equation}

In Figure \ref{fig:BB}, thimble $\mathcal{J}_{2n}$ is represented as a line connecting $z_{1}$ and $z_{2}$ in $w$ plane, where $w=x_{1}+ix_{2}$ and $\mathrm{Im} z_{1}>\mathrm{Im} z_{2}$.  After braided $\frac{\pi}{2}$ clockwise, thimble $\mathcal{J}_{2n}$ is twisted in the $2+1$ dimensional space time. It is shown in Figure \ref{fig:AB} after projected on $t,x_{2}$ plane. The knot invariant for this twist is $$<K>=\sum_{c,d}\mathcal{B}_{ab}^{cd}\mathcal{M}_{cd}.$$
Therefore, (\ref{braiding of B in NSB}) implies that
\begin{equation}
\sum_{c,d}\mathcal{B}_{ab}^{cd}\mathcal{M}_{cd}=q^{n}\mathcal{M}_{ab}.
\end{equation}

This means that the braiding of thimble $\mathcal{J}_{2n}$ without symmetry breaking  gives an eigenvalue of the braiding matrix in symmetry breaking case.

\begin{theorem}\label{Thmn2n}
For $\mathbbm{q}\neq 2n$,
\begin{equation}
\mathcal{B}\mathcal{J}_{f,\mathbbm{q}-f}=\left\{
                                             \begin{array}{ll}
                                               q^{- \frac{1}{2}}\mathcal{J}_{f,\mathbbm{q}-f}, & \hbox{$f=\mathbbm{q}-f$;} \\
                                               \mathcal{J}_{\mathbbm{q}-f,f}, & \hbox{$f > \mathbbm{q}-f$;} \\
                                               \mathcal{J}_{\mathbbm{q}-f,f}+(q^{-\frac{1}{2}}-q^{\frac{1}{2}})\mathcal{J}_{f,\mathbbm{q}-f}, & \hbox{$f< \mathbbm{q}-f$.}
                                             \end{array}
                                           \right.
\end{equation}
\end{theorem}

Proof:\begin{itemize}
        \item $f=\mathbbm{q}-f$: There is no wall-crossing. From Lemma \ref{lem:img part inv}, $$\mathcal{B}\mathcal{J}_{f,\mathbbm{q}-f}=q^{-\frac{1}{2}(\lambda^{f},\lambda^{f})}\mathcal{J}_{f,\mathbbm{q}-f}=q^{- \frac{1}{2}}\mathcal{J}_{f,\mathbbm{q}-f};$$
        \item  $f > \mathbbm{q}-f$: Also there is no wall-crossing. From Lemma \ref{lem:img part inv},
            $$\mathcal{B}\mathcal{J}_{f,\mathbbm{q}-f}=q^{-\frac{1}{2}(\lambda^{f},\lambda^{\mathbbm{q}-f})}\mathcal{J}_{f,\mathbbm{q}-f}=\mathcal{J}_{f,\mathbbm{q}-f};$$
        \item  $f < \mathbbm{q}-f$: From the conservation law of wall-crossing, there will be one wall-crossing term of $\mathcal{J}_{f,\mathbbm{q}-f}$. From Lemma \ref{lem:img part inv} and formula (\ref{braiding formula in SB}), $$\mathcal{B}\mathcal{J}_{\mathbbm{q}-f,f}=\mathcal{J}_{f,\mathbbm{q}-f},$$
            $$\mathcal{B}\mathcal{J}_{f,\mathbbm{q}-f}=q^{ -\frac{1}{2}(\lambda^{f},\lambda^{\mathbbm{q}-f})}\mathcal{J}_{\mathbbm{q}-f,f}+d\mathcal{J}_{f,\mathbbm{q}-f}=\mathcal{J}_{\mathbbm{q}-f,f}+d\mathcal{J}_{f,\mathbbm{q}-f},$$ where $d$ is an unknown constant. The transformation of $\mathcal{J}_{f,\mathbbm{q}-f}$ and $\mathcal{J}_{\mathbbm{q}-f,f}$ forms into a matrix:
          $$\mathcal{B}\left(
                         \begin{array}{c}
                           \mathcal{J}_{\mathbbm{q}-f,f} \\
                           \mathcal{J}_{f,\mathbbm{q}-f} \\
                         \end{array}
                       \right)=\left(
                                 \begin{array}{cc}
                                   0 & 1 \\
                                   1 & d \\
                                 \end{array}
                               \right)\left(
                                        \begin{array}{c}
                                          \mathcal{J}_{\mathbbm{q}-f,f} \\
                                          \mathcal{J}_{f,\mathbbm{q}-f} \\
                                        \end{array}
                                      \right).
$$
To determine $d$, we derive the braiding matrix of cycles $C_{\mathbbm{q}-f,f}$ and $C_{f,\mathbbm{q}-f}$. For convenience, we assume that $\mathbbm{q}-f=f+l$. From the second property of wall-crossing, the only possible transfer of simple roots is from $z_{2}$
 to $z_{1}$ , so the braiding of $C_{\mathbbm{q}-f,f}$ is easy: $$\mathcal{B}C_{\mathbbm{q}-f,f}=C_{f+l,f}=q^{-\frac{1}{2}(\lambda^{0},\lambda^{0})}C_{f,f+l}=q^{-\frac{1}{2}(\lambda^{0},\lambda^{0})}C_{f,\mathbbm{q}-f}=q^{-\frac{1}{2}}C_{f,\mathbbm{q}-f}.$$

\begin{figure}
\begin{center}
\includegraphics[width=7cm,height=5cm]{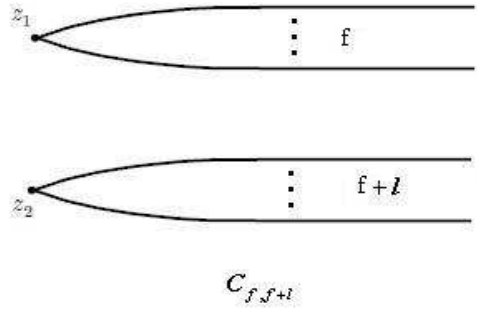}
\caption{\label{fig:Cff+l}$C_{f,f+l}$ before braiding}
\end{center}
\end{figure}

\begin{figure}
\begin{center}
\includegraphics[width=7.5cm,height=4cm]{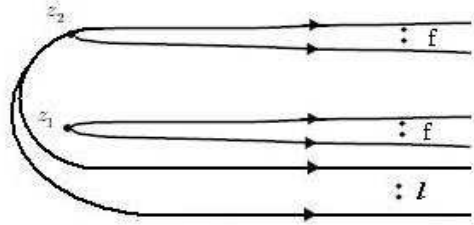}
\caption{\label{fig:WC1}$C_{f,f+l}$ after braiding}
\end{center}
\end{figure}

\begin{figure}
\begin{center}
\includegraphics[width=15cm,height=5cm]{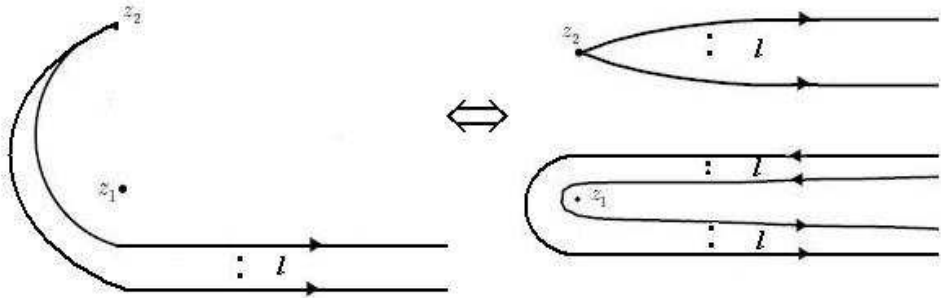}
\caption{\label{fig:WC2}Homology equivalence of wall-crossing part}
\end{center}
\end{figure}

The braiding of $C_{f,\mathbbm{q}-f}=C_{f,f+l}$ will cause wall-crossing. As is shown in Figure \ref{fig:Cff+l}, Figure \ref{fig:WC1} and Figure \ref{fig:WC2},  wall-crossing part is equivalent in homology to a zig-zag cycle, which is starting from  $z_{2}$, heads directly to $Re z = \infty$ before doubling back around $z_{1}$ and returning to $Re z = \infty$. Thus, there are three pieces in the wall-crossing part with two of them near $z_{1}$.

Consider the integral of the generalized Yang-Yang function on $C_{f,f+l}$:
\begin{equation}
\begin{split}
&\int _{C_{f,f+l}}e^{-\frac{\mathscr{W}}{k+h^{\vee }}}\prod _{r}dw_{r}\\
&=\int _{z_{1}}^{+\infty}...\int _{z_{1}}^{+\infty}\prod_{j=1}^{f}dw_{j}\int _{z_{2}}^{+\infty}...\int _{z_{2}}^{+\infty}\prod_{s=f+1}^{2f}dw_{s}\int _{z_{2}}^{+\infty}...\int _{z_{2}}^{+\infty}\prod_{t=2f+1}^{2f+l}dw_{t}\\
&((z_{1}-z_{2})^{\frac{(\lambda,\lambda)}{k+h^{\vee }}}\prod _{r=1}^{2f+l}\prod _{a=1}^{2}(w_{r}-z_{a})^{-\frac{(\alpha_{i_{r}},\lambda)}{k+h^{\vee }}}\prod _{r<r'}(w_{r}-w_{r'})^{\frac{(\alpha_{i_{r}},\alpha_{i_{r'}})}{k+h^{\vee }}}e^{\frac{c(\sum_{r}w_{r}-\frac{1}{2}\sum_{a}\parallel \lambda_{a}\parallel z_{a})}{k+h^{\vee }}})
\end{split}
\end{equation}

After braiding, this integral becomes:
\begin{equation}
\begin{split}
&\int _{z_{2}}^{+\infty}...\int _{z_{2}}^{+\infty}\prod_{j=1}^{f}dw_{j}\int _{z_{1}}^{+\infty}...\int _{z_{1}}^{+\infty}\prod_{s=f+1}^{2f}dw_{s}q^{-\frac{1}{2}[(\lambda,\lambda)-(\sum_{j}\alpha_{i_{j}},\lambda)-(\sum_{s}\alpha_{i_{s}},\lambda)+(\sum_{j}\alpha_{i_{j}},\sum_{s}\alpha_{i_{s}})]}\\
&(\int _{z_{1}}^{z_{2}}...\int _{z_{1}}^{z_{2}}\prod_{t=2f+1}^{2f+l}dw_{t}q^{-\frac{1}{2}[(\sum_{t}\alpha_{i_{t}},\sum_{j}\alpha_{i_{j}})+(\sum_{t}\alpha_{i_{t}},\sum_{s}\alpha_{i_{s}})-(\sum_{t}\alpha_{i_{t}},\lambda)-(\sum_{t}\alpha_{i_{t}},\lambda)]}+\int _{z_{2}}^{+\infty}...\int _{z_{2}}^{+\infty}\prod_{t=2f+1}^{2f+l}dw_{t})\\
&(z_{1}-z_{2})^{\frac{(\lambda,\lambda)}{k+h^{\vee }}}\prod _{r=1}^{2f+l}\prod _{a=1}^{2}(w_{r}-z_{a})^{-\frac{(\alpha_{i_{r}},\lambda)}{k+h^{\vee }}}\prod _{r<r'}(w_{r}-w_{r'})^{\frac{(\alpha_{i_{r}},\alpha_{i_{r'}})}{k+h^{\vee }}}e^{\frac{c(\sum_{r}w_{r}-\frac{1}{2}\sum_{a}\parallel \lambda_{a}\parallel z_{a})}{k+h^{\vee }}}\\
&=\int _{z_{2}}^{+\infty}...\int _{z_{2}}^{+\infty}\prod_{j=1}^{f}dw_{j}\int _{z_{1}}^{+\infty}...\int _{z_{1}}^{+\infty}\prod_{s=f+1}^{2f}dw_{s}q^{-\frac{1}{2}[(\lambda-\sum_{j}\alpha_{i_{j}},\lambda-\sum_{s}\alpha_{i_{s}})]}\\
&[(\int _{z_{1}}^{+\infty}...\int _{z_{1}}^{+\infty}-\int _{z_{2}}^{+\infty}...\int _{z_{2}}^{+\infty})\prod_{t=2f+1}^{2f+l}dw_{t}q^{\frac{1}{2}[(\sum_{t}\alpha_{i_{t}},\lambda-\sum_{j}\alpha_{i_{j}})+(\sum_{t}\alpha_{i_{t}},\lambda-\sum_{s}\alpha_{i_{s}})]}+\int _{z_{2}}^{+\infty}...\int _{z_{2}}^{+\infty}\prod_{t=2f+1}^{2f+l}dw_{t}]\\
&(z_{1}-z_{2})^{\frac{(\lambda,\lambda)}{k+h^{\vee }}}\prod _{r=1}^{2f+l}\prod _{a=1}^{2}(w_{r}-z_{a})^{-\frac{(\alpha_{i_{r}},\lambda)}{k+h^{\vee }}}\prod _{r<r'}(w_{r}-w_{r'})^{\frac{(\alpha_{i_{r}},\alpha_{i_{r'}})}{k+h^{\vee }}}e^{\frac{c(\sum_{r}w_{r}-\frac{1}{2}\sum_{a}\parallel \lambda_{a}\parallel z_{a})}{k+h^{\vee }}}\\
&=\int _{z_{2}}^{+\infty}...\int _{z_{2}}^{+\infty}\prod_{j=1}^{f}dw_{j}\int _{z_{1}}^{+\infty}...\int _{z_{1}}^{+\infty}\prod_{s=f+1}^{2f}dw_{s}q^{-\frac{1}{2}[(\lambda^{f},\lambda^{f})]}\\
&[(\int _{z_{1}}^{+\infty}...\int _{z_{1}}^{+\infty}-\int _{z_{2}}^{+\infty}...\int _{z_{2}}^{+\infty})\prod_{t=2f+1}^{2f+l}dw_{t}q^{\frac{1}{2}[2(\sum_{t}\alpha_{i_{t}},\lambda^{f})]}+\int _{z_{2}}^{+\infty}...\int _{z_{2}}^{+\infty}\prod_{t=2f+1}^{2f+l}dw_{t}]\\
&(z_{1}-z_{2})^{\frac{(\lambda,\lambda)}{k+h^{\vee }}}\prod _{r=1}^{2f+l}\prod _{a=1}^{2}(w_{r}-z_{a})^{-\frac{(\alpha_{i_{r}},\lambda)}{k+h^{\vee }}}\prod _{r<r'}(w_{r}-w_{r'})^{\frac{(\alpha_{i_{r}},\alpha_{i_{r'}})}{k+h^{\vee }}}e^{\frac{c(\sum_{r}w_{r}-\frac{1}{2}\sum_{a}\parallel \lambda_{a}\parallel z_{a})}{k+h^{\vee }}}\\
&=q^{-\frac{1}{2}[2(\lambda^{f+l},\lambda^{f})-(\lambda^{f},\lambda^{f})]}\int _{C_{f+l,f}}e^{-\frac{\mathscr{W}}{k+h^{\vee }}}\prod _{r}dw_{r}+(q^{-\frac{1}{2}(\lambda^{f},\lambda^{f})}-q^{-\frac{1}{2}[2(\lambda^{f+l},\lambda^{f})-(\lambda^{f},\lambda^{f})]})\int _{C_{f,f+l}}e^{-\frac{\mathscr{W}}{k+h^{\vee }}}\prod _{r}dw_{r}
\end{split}
\end{equation}
The ranges of indexes $a,j,r,r',s,t$ are $a=1,2$, $j=1,2,...,f$, $s=f+1,f+2,...,2f$, $t=2f+1,2f+2,...,2f+l$ and $r,r'=1,2,...2f+l$. In the first formula of equations above ,
$q^{-\frac{1}{2}[(\lambda,\lambda)-(\sum_{j}\alpha_{i_{j}},\lambda)-(\sum_{s}\alpha_{i_{s}},\lambda)+(\sum_{j}\alpha_{i_{j}},\sum_{s}\alpha_{i_{s}})]}$ comes from the braiding of the factor$(z_{1}-z_{2})^{\frac{(\lambda,\lambda)}{k+h^{\vee }}}$, $(w_{j}-z_{2})^{-\frac{(\alpha_{i_{j}},\lambda)}{k+h^{\vee }}}$, $(w_{s}-z_{1})^{-\frac{(\alpha_{i_{s}},\lambda)}{k+h^{\vee }}}$ and $(w_{j}-w_{s})^{\frac{(\alpha_{i_{j}},\alpha_{i_{s}})}{k+h^{\vee }}}$,  $q^{-\frac{1}{2}[(\sum_{t}\alpha_{i_{t}},\sum_{j}\alpha_{i_{j}})+(\sum_{t}\alpha_{i_{t}},\sum_{s}\alpha_{i_{s}})-(\sum_{t}\alpha_{i_{t}},\lambda)-(\sum_{t}\alpha_{i_{t}},\lambda)]}$ from $(w_{t}-w_{j})^{\frac{(\alpha_{i_{t}},\alpha_{i_{j}})}{k+h^{\vee }}}$,$(w_{t}-w_{s})^{\frac{(\alpha_{i_{t}},\alpha_{i_{s}})}{k+h^{\vee }}}$ , $(w_{t}-z_{2})^{-\frac{(\alpha_{i_{t}},\lambda)}{k+h^{\vee }}}$ and $(w_{t}-z_{1})^{-\frac{(\alpha_{i_{t}},\lambda)}{k+h^{\vee }}}$.

Thus, we have
\begin{equation}
\begin{split}
\mathcal{B}C_{f,\mathbbm{q}-f}&=\mathcal{B}C_{f,f+l}\\
                              &=q^{-\frac{1}{2}[2(\lambda^{f+l},\lambda^{f})-(\lambda^{f},\lambda^{f})]}(C_{f+l,f}-C_{f,f+l})+q^{-\frac{1}{2}(\lambda^{f},\lambda^{f})}C_{f,f+l}\\
                              &=q^{ \frac{1}{2}}(C_{f+l,f}-C_{f,f+l})+q^{ -\frac{1}{2}}C_{f,f+l}\\
                              &=q^{ \frac{1}{2}}C_{f+l,f}+(q^{ -\frac{1}{2}}-q^{ \frac{1}{2}})C_{f,f+l}\\
                              &=q^{ \frac{1}{2}}C_{\mathbbm{q}-f,f}+(q^{ -\frac{1}{2}}-q^{ \frac{1}{2}})C_{f,\mathbbm{q}-f},\\
\end{split}
\end{equation}

Thus,
$$\mathcal{B}\left(
                         \begin{array}{c}
                           C_{\mathbbm{q}-f,f} \\
                           C_{f,\mathbbm{q}-f} \\
                         \end{array}
                       \right)=\left(
                                 \begin{array}{cc}
                                   0 & q^{-\frac{1}{2}} \\
                                   q^{ \frac{1}{2}} & q^{ -\frac{1}{2}}-q^{ \frac{1}{2}} \\
                                 \end{array}
                               \right)\left(
                                        \begin{array}{c}
                                          C_{\mathbbm{q}-f,f} \\
                                          C_{f,\mathbbm{q}-f} \\
                                        \end{array}
                                      \right).
$$
$\{C_{\mathbbm{q}-f,f}, C_{f,\mathbbm{q}-f}\}$ and $\{\mathcal{J}_{\mathbbm{q}-f,f}, \mathcal{J}_{f,\mathbbm{q}-f}\}$ are two basis in the same vector space, so braiding matrices in these two basis are similar to each other. Thus $d=q^{ -\frac{1}{2}}-q^{ \frac{1}{2}}.$ This completes the proof.
      \end{itemize}

\begin{lemma}\label{Thm2n}
For $\mathbbm{q}=2n$,
\begin{equation}
\mathcal{B}\mathcal{J}_{f,2n-f}=\left\{
                                             \begin{array}{ll}
                                               \mathcal{J}_{2n-f,f}+\sum_{i=1}^{2n-f}\beta_{f,2n-f}^{2n-f-i,f+i}\mathcal{J}_{2n-f-i,f+i}, & \hbox{$f=2n-f$;} \\
                                               q^{ \frac{1}{2}}\mathcal{J}_{2n-f,f}+\sum_{i=1}^{2n-f}\beta_{f,2n-f}^{2n-f-i,f+i}\mathcal{J}_{2n-f-i,f+i}, & \hbox{$f \neq2n-f$,} \\
                                            \end{array}
                                           \right.
\end{equation}
\end{lemma}
where $\beta_{f,2n-f}^{2n-f-i,f+i}$ are unknown constants.

Proof:By the property of the wall-crossing , if we choose $\mathcal{J}_{2n,0}$, $\mathcal{J}_{2n-1,1}$ ,..., $\mathcal{J}_{0,2n}$ as a set of basis,  $\mathcal{B}$ will be a triangular matrix on $V_{2n}$. The skew diagonal elements are coming from the formula (\ref{braiding formula in SB}) of the braiding in the symmetry breaking case.

\begin{figure}
\begin{center}
\includegraphics[width=6cm,height=4cm]{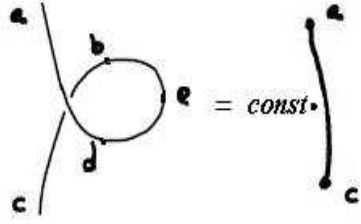}
\caption{\label{fig:Mcdt} Constraint for $\mathcal{M}$}
\end{center}
\end{figure}

To derive a knot invariant, we consider the general case of $\mathrm{dim} V_{\lambda}=m+1$ and assume that
for $\mathbbm{q}\neq m$,
\begin{equation}\label{equ:g braiding <m}
\mathcal{B}\mathcal{J}_{s,\mathbbm{q}-s}=\left\{
                                             \begin{array}{ll}
                                               \gamma\mathcal{J}_{s,\mathbbm{q}-s}, & \hbox{$s=\mathbbm{q}-s$;} \\
                                               \mathcal{J}_{\mathbbm{q}-s,s}, & \hbox{$s > \mathbbm{q}-s$;} \\
                                               \mathcal{J}_{\mathbbm{q}-s,s}+(\gamma-\gamma^{-1})\mathcal{J}_{s,\mathbbm{q}-s}, & \hbox{$s< \mathbbm{q}-s$,}
                                             \end{array}
                                           \right.
\end{equation}
and
\begin{equation}\label{equ:g braiding =m}
\mathcal{B}\mathcal{J}_{s,m-s}=\left\{
                                             \begin{array}{ll}
                                               \mathcal{J}_{m-s,s}+\sum_{i=1}^{2m-s}\beta_{s,m-s}^{m-s-i,s+i}\mathcal{J}_{m-s-i,s+i}, & \hbox{$s=m-s$;} \\
                                               \gamma^{-1}\mathcal{J}_{m-s,s}+\sum_{i=1}^{m-s}\beta_{s,m-s}^{m-s-i,s+i}\mathcal{J}_{m-s-i,s+i}, & \hbox{$s \neq m-s$.} \\
                                            \end{array}
                                           \right.
\end{equation}

 As is shown in Figure \ref{fig:Mcdt}, braiding and fusion matrices must satisfy the following condition:
\begin{equation}\label{equ:cond for M}
\sum_{b,d,e}\mathcal{B}_{cd}^{ab}\mathcal{M}_{be}\mathcal{M}^{de}=C\delta_{c}^{a}, \hbox{where $C$ is a constant}.
\end{equation}

We have:

\begin{lemma}\label{Lemma constraint}
\begin{description}
  \item[When $m$ is even] the condition (\ref{equ:cond for M}) implies that
\begin{equation}
\mathcal{M}_{am-a}\mathcal{M}^{am-a}=\left\{
                                              \begin{array}{ll}
                                                \gamma^{2a-m-1}, & \hbox{$\frac{m}{2}<a\leq m$;} \\
                                                1, & \hbox{$a=\frac{m}{2}$;} \\
                                                \gamma^{2a-m+1}, & \hbox{$0\leq a< \frac{m}{2}$,}
                                              \end{array}
                                            \right.
\end{equation} $$C=\gamma^{m}\hbox{ and  }\beta_{am-a}^{am-a}=(\gamma-\gamma^{-1})(1-\gamma^{2a-m+1}).$$
  \item[When $m$ is odd] the condition (\ref{equ:cond for M}) implies that
\begin{equation}
\mathcal{M}_{am-a}\mathcal{M}^{am-a}=\left\{
                                              \begin{array}{ll}
                                                 x\gamma^{2a-m-1}, & \hbox{$\frac{m}{2}<a\leq m$;} \\
                                                x^{-1}\gamma^{2a-m+1}, & \hbox{$0\leq a< \frac{m}{2}$,}
                                              \end{array}
                                            \right.
\end{equation} $$C=x\gamma^{m}\hbox{ and  }\beta_{am-a}^{am-a}=\gamma-\gamma^{-1}+(\gamma^{-1}-x^{-2}\gamma)\gamma^{2a-m+1},\quad 0\leq a<\frac{m}{2},$$ where $x=\mathcal{M}_{[\frac{m}{2}]+1[\frac{m}{2}]}\mathcal{M}^{[\frac{m}{2}]+1[\frac{m}{2}]}.$
\end{description}

\end{lemma}
Proof: The left hand side of (\ref{equ:cond for M}) is equal to
$$L.H.S=\sum_{e}\mathcal{B}_{cm-e}^{am-e}\mathcal{M}_{m-ee}\mathcal{M}^{m-ee}.$$ From the conservation law of wall-crossing,  $\mathcal{B}_{cd}^{ab}$ can be nonzero only when $c+d=a+b$. Thus,
$$L.H.S=\delta_{c}^{a}\sum_{e}\mathcal{B}_{am-e}^{am-e}\mathcal{M}_{m-ee}\mathcal{M}^{m-ee}.$$

\begin{itemize}
  \item When $\frac{m}{2}<a\leq m$: $$L.H.S=\delta_{c}^{a}\sum_{e\leq m-a}\mathcal{B}_{am-e}^{am-e}\mathcal{M}_{m-ee}\mathcal{M}^{m-ee}.$$ From (\ref{equ:g braiding <m}), $$L.H.S=\delta_{c}^{a}[\gamma\mathcal{M}_{am-a}\mathcal{M}^{am-a}+\sum_{e<m-a}(\gamma-\gamma^{-1})\mathcal{M}_{m-ee}\mathcal{M}^{m-ee}].$$
  \item When $a=\frac{m}{2}$: From (\ref{equ:g braiding <m}) and (\ref{equ:g braiding =m}), $$L.H.S=\delta_{c}^{a}[\mathcal{M}_{\frac{m}{2}\frac{m}{2}}\mathcal{M}^{\frac{m}{2}\frac{m}{2}}+\sum_{e<\frac{m}{2}}(\gamma-\gamma^{-1})\mathcal{M}_{m-ee}\mathcal{M}^{m-ee}].$$
  \item When $0\leq a<\frac{m}{2}$: From (\ref{equ:g braiding <m}) and (\ref{equ:g braiding =m}),
      \begin{equation}
      \begin{split}
       L.H.S&=\delta_{c}^{a}[\gamma\mathcal{M}_{am-a}\mathcal{M}^{am-a}+(\gamma-\gamma^{-1})\sum_{e<m-a,e\neq a}\mathcal{M}_{m-ee}\mathcal{M}^{m-ee}+\mathcal{B}_{am-a}^{am-a}\mathcal{M}_{m-aa}\mathcal{M}^{m-aa}]\\
&=\delta_{c}^{a}[\gamma\mathcal{M}_{am-a}\mathcal{M}^{am-a}+(\gamma-\gamma^{-1})\sum_{e<m-a,e\neq a}\mathcal{M}_{m-ee}\mathcal{M}^{m-ee}+\beta_{am-a}^{am-a}\mathcal{M}_{m-aa}\mathcal{M}^{m-aa}].\end{split}
\end{equation}
\end{itemize}

Thus,
\begin{itemize}
  \item When $\frac{m}{2}<a\leq m$: \begin{equation}\label{up}\gamma\mathcal{M}_{am-a}\mathcal{M}^{am-a}+\sum_{e<m-a}(\gamma-\gamma^{-1})\mathcal{M}_{m-ee}\mathcal{M}^{m-ee}=C.\end{equation}
  \item When $a=\frac{m}{2}$: \begin{equation}\label{mid}\mathcal{M}_{\frac{m}{2}\frac{m}{2}}\mathcal{M}^{\frac{m}{2}\frac{m}{2}}+\sum_{e<\frac{m}{2}}(\gamma-\gamma^{-1})\mathcal{M}_{m-ee}\mathcal{M}^{m-ee}=C.\end{equation}
  \item When $0\leq a<\frac{m}{2}$: \begin{equation}\label{dwn}\gamma\mathcal{M}_{am-a}\mathcal{M}^{am-a}+(\gamma-\gamma^{-1})\sum_{e<m-a,e\neq a}\mathcal{M}_{m-ee}\mathcal{M}^{m-ee}+\beta_{am-a}^{am-a}\mathcal{M}_{m-aa}\mathcal{M}^{m-aa}=C.\end{equation}
\end{itemize}
\begin{description}
  \item[When $m$ is even] From (\ref{up}) and (\ref{mid}),
$$\frac{\mathcal{M}_{am-a}\mathcal{M}^{am-a}}{\mathcal{M}_{a-1m-a+1}\mathcal{M}^{a-1m-a+1}}=\left\{
                                                                    \begin{array}{ll}
                                                                      \gamma ^{2}, & \hbox{$\frac{m}{2}+2\leq a\leq m$;} \\
                                                                     \gamma , & \hbox{$a=\frac{m}{2}+1$.}
                                                                    \end{array}
                                                                  \right.
$$
It should be noted that $$\mathcal{M}_{am-a}\mathcal{M}^{am-a}=(\mathcal{M}_{m-a a}\mathcal{M}^{m-a a})^{-1}.$$ Especially, $$\mathcal{M}_{\frac{m}{2}\frac{m}{2}}\mathcal{M}^{\frac{m}{2}\frac{m}{2}}=1.$$
This leads to $$\mathcal{M}_{am-a}\mathcal{M}^{am-a}=\left\{
                                              \begin{array}{ll}
                                                \gamma^{2a-m-1}, & \hbox{$\frac{m}{2}<a\leq m$;} \\
                                                1, & \hbox{$a=\frac{m}{2}$;} \\
                                                \gamma^{2a-m+1}, & \hbox{$0\leq a< \frac{m}{2}$.}
                                              \end{array}
                                            \right.$$
Let $a=m$, from (\ref{up}), we have $$C=\gamma\mathcal{M}_{m0}\mathcal{M}^{m0}=\gamma^{m}.$$
Since all of $\mathcal{M}_{m-aa}\mathcal{M}^{m-aa}$ are known, from (\ref{dwn}), $$\beta_{am-a}^{am-a}=(\gamma-\gamma^{-1})(1-\gamma^{2a-m+1}).$$
  \item[When $m$ is odd]  From (\ref{up}),
$$\frac{\mathcal{M}_{am-a}\mathcal{M}^{am-a}}{\mathcal{M}_{a-1m-a+1}\mathcal{M}^{a-1m-a+1}}=\gamma ^{2},\quad a-1>\frac{m}{2}.
$$
It should be noted that $$\mathcal{M}_{am-a}\mathcal{M}^{am-a}=(\mathcal{M}_{m-a a}\mathcal{M}^{m-a a})^{-1}.$$
If we assume that $$x=\mathcal{M}_{[\frac{m}{2}]+1[\frac{m}{2}]}\mathcal{M}^{[\frac{m}{2}]+1[\frac{m}{2}]},$$ then $$\mathcal{M}_{[\frac{m}{2}][\frac{m}{2}]+1}\mathcal{M}^{[\frac{m}{2}][\frac{m}{2}]+1}=x^{-1}.$$
This leads to $$\mathcal{M}_{am-a}\mathcal{M}^{am-a}=\left\{
                                              \begin{array}{ll}
                                                 x\gamma^{2a-m-1}, & \hbox{$\frac{m}{2}<a\leq m$;} \\
                                                x^{-1}\gamma^{2a-m+1}, & \hbox{$0\leq a< \frac{m}{2}$.}
                                              \end{array}
                                            \right.$$
Let $a=m$, from (\ref{up}), we have $$C=\gamma\mathcal{M}_{m0}\mathcal{M}^{m0}=x\gamma^{m}.$$
Since all of $\mathcal{M}_{m-aa}\mathcal{M}^{m-aa}$ are known, from (\ref{dwn}), $$\beta_{am-a}^{am-a}=\gamma-\gamma^{-1}+(\gamma^{-1}-x^{-2}\gamma)\gamma^{2a-m+1},\quad 0\leq a<\frac{m}{2}.$$
\end{description}

This completes the proof.

For the fundamental representation of $B_{n}$ Lie algebra, $m=2n$, $\gamma=q^{-\frac{1}{2}}$.  From Lemma \ref{Lemma constraint}, \begin{equation}
\mathcal{M}_{a2n-a}\mathcal{M}^{a2n-a}=\left\{
                                              \begin{array}{ll}
                                               q^{n-a+\frac{1}{2}}, & \hbox{$n<a\leq 2n$;} \\
                                                1, & \hbox{$a=n$;} \\
                                               q^{n-a-\frac{1}{2}}, & \hbox{$0\leq a< n$,}
                                              \end{array}
                                            \right.
\end{equation} $$C=q^{-n}\hbox{ and  }\beta_{a2n-a}^{a2n-a}=(q^{-\frac{1}{2}}-q^{\frac{1}{2}})(1-q^{n-a-\frac{1}{2}}).$$

\begin{lemma}\label{Lem tr det}
For $\mathbbm{q}=2n$, if we choose $\mathcal{J}_{2n,0}$, $\mathcal{J}_{2n-1,1}$, ..., $\mathcal{J}_{1,2n-1}$, $\mathcal{J}_{0,2n}$ as a set of basis, then $\mathcal{B}$ is a $(2n+1)\times(2n+1)$ triangular matrix on $V_{2n}$,
\begin{equation}
\mathrm{Tr}\mathcal{B}=n(q^{-\frac{1}{2}}-q^{\frac{1}{2}})+q^{n}
\end{equation}
and
\begin{equation}
\mathrm{det}\mathcal{B}=(-1)^{n}q^{n}.
\end{equation}
\end{lemma}
Proof: From Lemma \ref{Thm2n}, the skew diagonal elements are known, therefore
$$\mathrm{det}\mathcal{B}=(-1)^{n}\prod_{0\leq a\leq 2n;a\neq n} q^{\frac{1}{2}}=(-1)^{n}q^{n}.$$
From Lemma \ref{Thm2n} and Lemma \ref{Lemma constraint},
$$\mathrm{Tr}\mathcal{B}=1+\sum_{0\leq a\leq n-1}\beta_{a2n-a}^{a2n-a}=n(q^{-\frac{1}{2}}-q^{\frac{1}{2}})+q^{n}.$$
This completes the proof.

\begin{figure}
\begin{center}
\includegraphics[width=8cm,height=4cm]{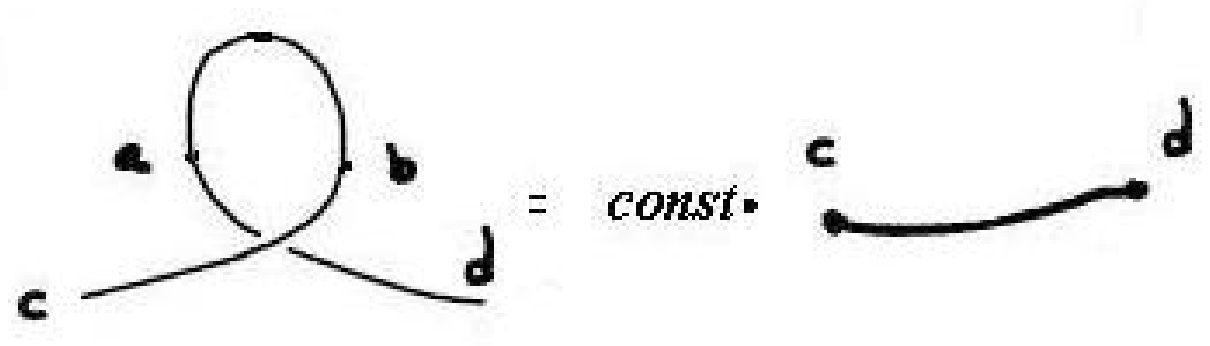}
\caption{\label{fig:Mcdt2}}
\end{center}
\end{figure}

After rotating Figure \ref{fig:Mcdt} $\frac{\pi}{2}$ counter clockwise, we have Figure \ref{fig:Mcdt2}, i.e. \begin{equation}\label{equ: fusion matrix as eigenvector of inverse}
\sum_{a,b}(\mathcal{B}^{-1})_{cd}^{ab}\mathcal{M}_{ab}=C\mathcal{M}_{cd},\end{equation} where $C$ must be the same constant in Figure \ref{fig:Mcdt}.
Thus, \begin{equation}\label{equ: fusion matrix as eigenvector}
\sum_{a,b}\mathcal{B}_{cd}^{ ab}\mathcal{M}_{ab}=C^{-1}\mathcal{M}_{cd}.
\end{equation}
Since $$\mathcal{M}_{ab}=\delta^{b2n-a}\mathcal{M}_{a2n-a},$$ if we denote $$\zeta^{t}=(
                                                                          \begin{array}{ccccc}
                                                                            \zeta_{0} & \zeta_{1} & \zeta_{2} & ... & \zeta_{2n} \\
                                                                          \end{array}
                                                                        ), \zeta_{a}=\mathcal{M}_{a2n-a},$$
then \begin{equation}\label{equ: fusion as eigenvector}
\mathcal{B}\zeta=C^{-1}\zeta\hbox{ on }V_{2n},
\end{equation} i.e. $\zeta$ is an eigenvector of $\mathcal{B}$ on $V_{2n}$ with respect to the eigenvalue $C^{-1}=q^{n}$.
This is in accordance with the braiding formula (\ref{braiding of B in NSB}) for thimble $\mathcal{J}_{2n}$.

To find other eigenvalues of $\mathcal{B}$ on $V_{2n}$, we introduce the following Lemma:

\begin{lemma}\label{prop Laurent ring }

Let $a_{i}\in$ $\mathbb{Z}[y,y^{-1}]$, $\mathbb{Z}[y,y^{-1}]$ is a ring of Laurent polynomials of $y$. If

\begin{equation}
\left\{
                              \begin{array}{ll}
                                a_{1}\cdot a_{2}\cdot a_{3}\cdot...a_{m}=(-1)^{m_{2}}; \\
                                a_{1}+ a_{2}+... a_{m}=m_{1}\cdot y-m_{2}\cdot y^{-1},
                              \end{array}
                            \right.
\end{equation}where $m_{1}, m_{2}\in \mathbb{Z}_{+}$ and $m_{1}+ m_{2}=m$,
then \begin{equation}
a_{i}=\left\{
                                                  \begin{array}{ll}
                                                    y, & \hbox{$i=1,2,...,m_{1}$;} \\
                                                    -y^{-1}, & \hbox{$i=m_{1}+1,m_{1}+2,...,m,$}
                                                  \end{array}
                                                \right.
\end{equation}
 up to a symmetry group $S_{m}$ action.
\end{lemma}
Proof:If for some $i$, $a_{i}$ is not a monomial in $\mathbb{Z}[y,y^{-1}]$, then from $$a_{1}\cdot a_{2}\cdot a_{3}\cdot... a_{m}=(-1)^{m_{2}},$$  we have $$a_{1}\cdot a_{2}\cdot a_{3}\cdot...\breve{a_{i}}\cdot... a_{m}=\frac{(-1)^{m_{2}}}{a_{i}}\notin \mathbb{Z}[y,y^{-1}].$$ But for $j\ne i$, $a_{j}\in \mathbb{Z}[y,y^{-1}]$, so $$a_{1}\cdot a_{2}\cdot a_{3}\cdot...\breve{a_{i}}\cdot... a_{m}\in \mathbb{Z}[y,y^{-1}].$$
The assumption leads to a contradiction, thus for any $i=1,2,...m$, $a_{i}$ is a monomial in $\mathbb{Z}[y,y^{-1}]$.
Assume that $$a_{i}=d_{i}y^{l_{i}},\quad d_{i}\in \mathbb{Z}, l_{i}\in \mathbb{Z}.$$ Then from $$a_{1}\cdot a_{2}\cdot a_{3}\cdot... a_{m}=(-1)^{m_{2}},$$ we know \begin{equation}
\left\{
         \begin{array}{ll}
           d_{i}=\pm1 \\
           \sum_{i}l_{i}=0.
         \end{array}
       \right.
\end{equation}
Let $$S_{+}=\{i|l_{i}=1\}$$ and $$S_{-}=\{i|l_{i}=-1\},$$ then
\begin{equation}
\left\{
         \begin{array}{ll}
          \sum_{i\in S_{+}} d_{i}=m_{1} \\
          \sum_{i\in S_{-}} d_{i}=-m_{2}.
         \end{array}
       \right.
\end{equation}
From $d_{i}=\pm1$, we have $$\#S_{+}\geqslant m_{1}\hbox{ and }\#S_{-}\geqslant m_{2}.$$
$S_{+}$ and $S_{-}$ are subsets of $\{1,2,3...,m\}$, therefore $$\#S_{+}+\#S_{-}\leqslant m.$$
Thus, $$\#S_{+}= m_{1}\hbox{ and }\#S_{-}= m_{2};$$
$$d_{i}=\left\{
           \begin{array}{ll}
             1, & \hbox{$i\in S_{+}$;} \\
             -1, & \hbox{$i\in S_{-}$.}
           \end{array}
         \right.
$$
This completes the proof.

\begin{theorem}
All $2n+1$ eigenvalues for $\mathcal{B}$ on $V_{2n}$ are $$e_{i}=\left\{
                                                                   \begin{array}{ll}
                                                                     q^{-\frac{1}{2}}, & \hbox{$1\leq i\leq n$;} \\
                                                                     -q^{\frac{1}{2}}, & \hbox{$n+1\leq i\leq 2n$;} \\
                                                                     q^{n}, & \hbox{$i=2n+1$.}
                                                                   \end{array}
                                                                 \right.$$
\end{theorem}
Proof: We know the trace and determinant of $\mathcal{B}$ on $V_{2n}$ from Lemma \ref{Lem tr det}:
$$\mathrm{Tr}\mathcal{B}=n(q^{-\frac{1}{2}}-q^{\frac{1}{2}})+q^{n},$$
$$\mathrm{det}\mathcal{B}=(-1)^{n}q^{n}.$$
From Figure \ref{fig:Mcdt2}, one of $2n+1$ eigenvalues is $e_{2n+1}=q^{n}$. Thus, $\sum_{i=1}^{2n}e_{i}=n\cdot(q^{-\frac{1}{2}}-q^{\frac{1}{2}})$ and $\prod_{i=1}^{2n}e_{i}=(-1)^{n}$. From Lemma \ref{prop Laurent ring }, theorem is straight forward.

Now we introduce a lemma to derive the braiding matrix.
\begin{lemma}\label{lem:linear alg}
Let $S\in M(l, \mathbb{Z}[y,y^{-1}])$ be any $l\times l$ symmetric matrix of the shape \begin{equation}
S_{i,j}=\left\{
                    \begin{array}{ll}
                      0, & \hbox{$i<l+1-j$;}\\
                      a^{-1}, & \hbox{$i\neq j \quad\mathrm{and}\quad i=l+1-j$;} \\
                      1, & \hbox{$i=j \quad\mathrm{and}\quad i=l+1-j$;} \\
                      nonzero, & \hbox{$i>l+1-j$.}
                    \end{array}
                  \right.
\end{equation} Assume that it has only three different monomial eigenvalues: $a$ of multiplicity $m$, $-a^{-1}$ of multiplicity $l-m-1$, $c$ of multiplicity $1$. $a,-a^{-1},c,\frac{c-c^{-1}}{a-a^{-1}}\in \mathbb{Z}[y,y^{-1}]$. $\xi^{t}=(
                                                                          \begin{array}{cccc}
                                                                            \xi_{1} & \xi_{2} & ... & \xi_{l} \\
                                                                          \end{array}
                                                                        )
 $  is an eigenvector of $S$ satisfying
\begin{equation}S\cdot\xi=c\xi,\end{equation} and
\begin{equation}\label{norm}|\xi|^{2}=1-\frac{c-c^{-1}}{a-a^{-1}}.\end{equation}
Then
\begin{enumerate}
  \item $$S-S^{-1}=(a^{-1}-a)(\xi\cdot\xi^{t}-I).$$
  \item $$\xi_{i}\xi_{l+1-i}=1,\quad i=1,2,...,l
$$
\end{enumerate}
\end{lemma}
Proof:  Notice that all the conditions and results are in $\mathbb{Z}[y,y^{-1}]$. However, we prove this lemma in the field of fractions of the polynomial ring $\mathbb{Z}[y,y^{-1}]$. First, there is an orthogonal matrix $T$, $S=T\cdot\Lambda \cdot T^{t}$, where
\begin{equation}
\Lambda=\left(
     \begin{array}{ccccccc}
       a &  &  &  &  &  &  \\
        & ... &  &  &  &  &  \\
        &  & a &  &  &  &  \\
        &  &  & -a^{-1} &  &  &  \\
        &  &  &  & ... &  &  \\
        &  &  &  &  & -a^{-1} &  \\
        &  &  &  &  &  & c \\
     \end{array}
   \right).
\end{equation}
Obviously, \begin{equation}
\Lambda-\Lambda^{-1}=(a^{-1}-a)(\eta\cdot\eta^{t}-I), \quad\mathrm{where}\quad
\eta^{t}=(
             \begin{array}{cccc}
               0 & ... & 0 & (1-\frac{c-c^{-1}}{a-a^{-1}})^{\frac{1}{2}} \\
             \end{array}
           )
.
\end{equation}
Thus,
\begin{equation}\label{decomposition}
S-S^{-1}=(a^{-1}-a)(\xi'\cdot\xi'^{t}-I),\quad\mathrm{where}\quad\xi'=T\eta.
\end{equation}

\begin{equation}\label{xi is eigenvector}
 S\cdot\xi'=T\cdot\Lambda \cdot T^{t}\cdot T\cdot\eta=T\cdot\Lambda\cdot\eta=cT\cdot\eta=c\xi',
\end{equation}
i.e. $\xi'$ is an eigenvector of $S$ with respect to the eigenvalue $c$.

\begin{equation}
(S-S^{-1})\cdot\xi'=(c-c^{-1})\xi'=(a^{-1}-a)(\xi'\cdot\xi'^{t}\cdot\xi'-\xi')=(a^{-1}-a)(|\xi'|^{2}-1)\xi'.
\end{equation}
Thus, \begin{equation}
|\xi'|^{2}=1-\frac{c-c^{-1}}{a-a^{-1}}.
\end{equation}

Because the multiplicity of eigenvalue $c$ is $1$, the dimension of eigenvector space with respect to $c$ is also $1$. $\xi$ and $\xi'$ are linearly relative. $$|\xi'|^{2}=|\xi|^{2}=1-\frac{c-c^{-1}}{a-a^{-1}},$$ therefore $$\xi'=\pm\xi.$$

Thus, $$S-S^{-1}=(a^{-1}-a)(\xi\cdot\xi^{t}-I).$$

From
\begin{equation}
(S^{-1})_{i,j}=\left\{
                    \begin{array}{ll}
                      0, & \hbox{$i>l+1-j$;}\\
                      a, & \hbox{$i\neq j \quad\mathrm{and}\quad i=l+1-j$;} \\
                      1, & \hbox{$i=j \quad\mathrm{and}\quad i=l+1-j$;} \\
                      nonzero, & \hbox{$i<l+1-j$,}
                    \end{array}
                  \right.
\end{equation} it is straightforward to see \begin{equation}\xi\cdot\xi^{t}=\frac{1}{a^{-1}-a}(S-S^{-1})+I=\left(
                      \begin{array}{ccccc}
                         &  &  &  & 1 \\
                         & * &  & 1 &  \\
                         &  & ... &  &  \\
                         & 1 &  & * &  \\
                        1 &  &  &  &  \\
                      \end{array}
                    \right).
\end{equation}

Thus,
\begin{equation}\label{inverse}
\xi_{i}\xi_{l+1-i}=1.
\end{equation}
This completes the proof.

\begin{remark}
Lemma \ref{lem:linear alg} is used to prove that the braiding and fusion matrices we derived give Kauffman polynomial. The condition (\ref{norm}) is actually a constraint condition for the fusion matrix $\mathcal{M}$:$\sum_{a,b}(\mathcal{M}_{ab}\mathcal{M}^{ab})=1-\frac{c-c^{-1}}{a-a^{-1}}$. We will see it is equivalent to the condition $D_{\bigcirc }=((\alpha-\alpha^{-1})/z)+1$ in the definition of Kauffman polynomial.
\end{remark}

\begin{theorem}\label{B}
The braiding matrix $\mathcal{B}$ and fusion matrix $\mathcal{M}$ corresponding to the fundamental representation of $B_{n}$ Lie algebra satisfy the following conditions:
\begin{enumerate}
  \item $\mathcal{B}_{ac}^{bd}-(\mathcal{B}^{-1})_{ac}^{bd}=(q^{\frac{1}{2}}-q^{-\frac{1}{2}})(\mathcal{M}_{ac}\mathcal{M}^{bd}-\delta_{a}^{b}\delta_{c}^{d}).$
  \item $\sum_{a,b}(\mathcal{M}_{ab}\mathcal{M}^{ab})=\frac{q^{n}-q^{-n}}{q^{\frac{1}{2}}-q^{-\frac{1}{2}}}+1.$
  \item $\sum_{c,d}\mathcal{B}^{cd}_{ab}\mathcal{M}_{cd}=q^{n}\mathcal{M}_{ab}.$
  \item $\sum_{c,d}(\mathcal{B}^{-1})^{cd}_{ab}\mathcal{M}_{cd}=q^{-n}\mathcal{M}_{ab}.$
\end{enumerate}
\end{theorem}
Proof:
Now let $y=q^{\frac{1}{2}}$, $a=q^{-\frac{1}{2}}$, $c=q^{n}$, $l=2n+1$ and $S=\mathcal{B}$ on $V_{2n}$, then from Lemma \ref{lem:linear alg},  if
$$\mathcal{B}\cdot\xi=q^{n}\xi\hbox{   on }V_{2n},$$ and
$$|\xi|^{2}=1-\frac{q^{n}-q^{-n}}{q^{-\frac{1}{2}}-q^{\frac{1}{2}}},$$
then
  \begin{equation}\label{skein of xi}\mathcal{B}-\mathcal{B}^{-1}=(q^{\frac{1}{2}}-q^{-\frac{1}{2}})(\xi\cdot\xi^{t}-I)\hbox{   on }V_{2n}.\end{equation}
From (\ref{equ: fusion as eigenvector}), $$\zeta_{a}=\mathcal{M}_{a2n-a}$$ is also an eigenvector of $\mathcal{B}$ on $V_{2n}$ with respect to the eigenvalue $q^{n}$, thus $$\zeta=d\cdot\xi,$$ i.e. $$\mathcal{M}_{a2n-a}=d\cdot\xi_{a+1},$$ where $d$ is a constant.
From Lemma \ref{lem:linear alg}, $$\mathcal{M}^{2n-a a}=d^{-1}\cdot(\xi_{a+1})^{-1}=d^{-1}\cdot\xi_{2n+1-a},$$ or
$$\mathcal{M}^{a 2n-a }=d^{-1}\cdot\xi_{a+1}.$$
Thus $$\mathcal{M}_{a 2n-a }\mathcal{M}^{b 2n-b}=\xi_{a+1}\cdot\xi_{b+1},$$
$$\sum_{a,b}(\mathcal{M}_{ab}\mathcal{M}^{ab})=|\xi|^{2}=1-\frac{c-c^{-1}}{a-a^{-1}}=1+\frac{q^{n}-q^{-n}}{q^{\frac{1}{2}}-q^{-\frac{1}{2}}}.$$
From (\ref{skein of xi}),
 \begin{equation}\label{skein of B on V2n}
\mathcal{B}_{a2n-a}^{b2n-b}-(\mathcal{B}^{-1})_{a2n-a}^{b2n-b}=(q^{\frac{1}{2}}-q^{-\frac{1}{2}})(\mathcal{M}_{a 2n-a }\mathcal{M}^{b2n-b}-\delta_{a}^{b}\delta_{2n-a}^{2n-b})\hbox{   on }V_{2n}.\end{equation}

 For $\mathbbm{q}\neq 2n$,
\begin{equation}
\mathcal{B}\mathcal{J}_{m,\mathbbm{q}-m}=\left\{
                                             \begin{array}{ll}
                                               q^{- \frac{1}{2}}\mathcal{J}_{m,\mathbbm{q}-m}, & \hbox{$m=\mathbbm{q}-m$;} \\
                                               \mathcal{J}_{\mathbbm{q}-m,m}, & \hbox{$m > \mathbbm{q}-m$;} \\
                                               \mathcal{J}_{\mathbbm{q}-m,m}+(q^{-\frac{1}{2}}-q^{\frac{1}{2}})\mathcal{J}_{m,\mathbbm{q}-m}, & \hbox{$m< \mathbbm{q}-m$.}
                                             \end{array}
                                           \right.
\end{equation}
\begin{equation}
\mathcal{B}^{-1}\mathcal{J}_{m,\mathbbm{q}-m}=\left\{
                                             \begin{array}{ll}
                                               q^{ \frac{1}{2}}\mathcal{J}_{m,\mathbbm{q}-m}, & \hbox{$m=\mathbbm{q}-m$;} \\
                                               \mathcal{J}_{\mathbbm{q}-m,m}+(q^{\frac{1}{2}}-q^{-\frac{1}{2}})\mathcal{J}_{m,\mathbbm{q}-m}, & \hbox{$m > \mathbbm{q}-m$;} \\
                                               \mathcal{J}_{\mathbbm{q}-m,m}, & \hbox{$m< \mathbbm{q}-m$.}
                                             \end{array}
                                           \right.
\end{equation}
Thus, for $\mathbbm{q}\neq 2n$,
\begin{equation}\label{Bn1}
\mathcal{B}-\mathcal{B}^{-1}=(q^{-\frac{1}{2}}-q^{\frac{1}{2}})I\quad on\quad V_{\mathbbm{q}}.
\end{equation}

Combining two equations (\ref{Bn1}) and (\ref{skein of B on V2n}) together, we have
\begin{equation}
\mathcal{B}_{ac}^{bd}-(\mathcal{B}^{-1})_{ac}^{bd}=(q^{\frac{1}{2}}-q^{-\frac{1}{2}})(\mathcal{M}_{ac} \mathcal{M}^{bd}-\delta_{a}^{b}\delta_{c}^{d}).
\end{equation}

Since $C=q^{-n}$, (\ref{equ: fusion matrix as eigenvector} ) and (\ref{equ: fusion matrix as eigenvector of inverse} )
are equivalent to
\begin{equation}\label{}
  \sum_{c,d}\mathcal{B}^{cd}_{ab}\mathcal{M}_{cd}=q^{n}\mathcal{M}_{ab}
\end{equation}
and
\begin{equation}\label{}
  \sum_{c,d}(\mathcal{B}^{-1})^{cd}_{ab}\mathcal{M}_{cd}=q^{-n}\mathcal{M}_{ab}.
\end{equation}
This completes the proof.

To derive the braiding matrix $\mathcal{B}$ on $V_{2n}$, we choose a special solution of $\mathcal{M}_{ab}$. As the same as in the case of $A_{n}$ Lie algebra \cite{HL}, we choose $\mathcal{M}_{ab}$ satisfying the following normalization condition:
\begin{equation}\label{Norm for B ele}
\mathcal{M}_{a2n-a}\cdot\mathcal{M}_{2n-aa}=1,\quad a=0,1,2,...,2n,
\end{equation}
or equivalently,
\begin{equation}\label{Norm for B}
\mathcal{M}_{ab}=\mathcal{M}^{ab},\quad a,b=0,1,2,...,2n.
\end{equation}

It should be noticed that Theorem \ref{B} is independent of this normalization.

\begin{lemma}\label{Lemma constraint after norm}
From the condition (\ref{equ:cond for M}) $$\sum_{b,d,e}\mathcal{B}_{cd}^{ab}\mathcal{M}_{be}\mathcal{M}^{de}=C\delta_{c}^{a}$$ and (\ref{Norm for B}) $$\mathcal{M}_{ab}=\mathcal{M}^{ab},$$  we have
\begin{equation}
(\mathcal{M}_{a2n-a})^{2}=\left\{
                                              \begin{array}{ll}
                                                q^{\frac{2n-2a+1}{2}}, & \hbox{$n+1\leq a\leq 2n$;} \\
                                                1, & \hbox{$a=n$;} \\
                                                q^{\frac{2n-2a-1}{2}}, & \hbox{$0\leq a\leq n-1$,}
                                              \end{array}
                                            \right.
\end{equation}
\end{lemma}
Proof: It is straight forward from Lemma \ref{Lemma constraint}.

For convenience, we choose
\begin{equation}\label{equ:fusion B}
\mathcal{M}\left(
                    \begin{array}{c}
                      v_{2n} \\
                      v_{2n-1} \\
                      ... \\
                      v_{n+1}\\
                      v_{n} \\
                      v_{n-1} \\
                      ... \\
                      v_{1} \\
                      v_{0} \\
                    \end{array}
                  \right)
=\left(
                     \begin{array}{ccccccccc}
                        &  &  &  &   &   &   &   & q^{\frac{1}{2}(-n+\frac{1}{2})} \\
                        &  &  &  &  &  &  & q^{\frac{1}{2}(-n+1+\frac{1}{2})} &  \\
                        &  & 0&  &  &  & ... &  &  \\
                        &  &  &  &  & q^{-\frac{1}{4}} &  &  &  \\
                        &  &  &  & 1 &  &  &  &  \\
                        &  &  & q^{\frac{1}{4}} &  &  &  &  &  \\
                        &  & ... &  &  &  & 0 &  &  \\
                        & q^{\frac{1}{2}(n-1-\frac{1}{2})} &  &  &  &  &  &  &  \\
                       q^{\frac{1}{2}(n-\frac{1}{2})} &  &  &  &  &  &  &  &  \\
                     \end{array}
                   \right)
\left(
                    \begin{array}{c}
                      v_{2n} \\
                      v_{2n-1} \\
                      ... \\
                      v_{n+1}\\
                      v_{n} \\
                      v_{n-1} \\
                      ... \\
                      v_{1} \\
                      v_{0} \\
                    \end{array}
                  \right).
\end{equation}

Now we can use the formula (\ref{skein of B on V2n}) to determine the unknown elements $\beta_{a2n-a}^{2n-b b}$ in the braiding matrix $\mathcal{B}$ on $V_{2n}$.

\begin{theorem}For any $a=0,1,2,...2n-1$, $a<b$,
\begin{equation}
\beta_{a2n-a}^{2n-b b}=\left\{
                              \begin{array}{ll}
                                (q^{\frac{1}{2}}-q^{-\frac{1}{2}})q^{-\frac{1}{2}(a-b)}, & \hbox{$(a-n)(b-n)>0$;} \\
                                (q^{\frac{1}{2}}-q^{-\frac{1}{2}})q^{-\frac{1}{2}(a-b+\frac{1}{2})}, & \hbox{$(a-n)(b-n)=0$;} \\
                                (q^{\frac{1}{2}}-q^{-\frac{1}{2}})q^{-\frac{1}{2}(a-b+1)}+\delta_{a}^{2n-b}(q^{-\frac{1}{2}}-q^{\frac{1}{2}}), & \hbox{$(a-n)(b-n)<0$.}
                              \end{array}
                            \right.
\end{equation}
\end{theorem}
Proof: It is a straight forward calculation from (\ref{skein of B on V2n}) and (\ref{equ:fusion B}).

Thus, the braiding matrix $\mathcal{B}$ is derived and it is independent of the normalization condition (\ref{Norm for B ele}). One can check that when $e^{h}=q^{-\frac{1}{2}}$, this result is the same as that in section 7.3C of \cite{C}, therefore, braiding matrices we derived satisfy Yang-Baxter equation.

\subsection{Fundamental representation of $C_{n}$}

\begin{figure}
\begin{center}
\includegraphics[width=2cm,height=10cm]{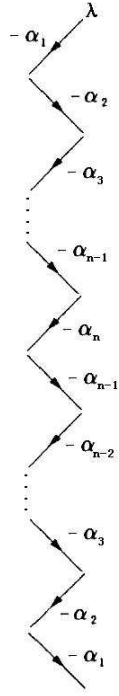}
\caption{\label{fig:Cn} Weights for the fundamental representation of $C_{n}$}
\end{center}
\end{figure}

$C_{n}$ has $n$ simple roots $\alpha_{i}, i=1,2,...,n$. Cartan matrix of $C_{n}$ is$$\left(
                                                     \begin{array}{ccccc}
                                                       2 & -1 &   &   &   \\
                                                       -1 & 2 & -1 &   &   \\
                                                         & -1 & ... &   &   \\
                                                         &   &   & 2 & -1 \\
                                                         &   &   & -2 & 2 \\
                                                     \end{array}
                                                   \right).$$
The highest weight of the fundamental representation $V_{\lambda}$ of $C_{n}$ is  $\lambda=(1,0,...,0)$. Weights of the fundamental representation are $\lambda,\lambda-\alpha_{1},\lambda-\alpha_{1}-\alpha_{2},...,\lambda-\sum_{i=1}^{n}\alpha_{i}, \lambda-\sum_{i=1}^{n}\alpha_{i}-\alpha_{n-1},\lambda-\sum_{i=1}^{n}\alpha_{i}-\alpha_{n-1}-\alpha_{n-2},...,\lambda-\sum_{i=1}^{n-1}2\alpha_{i}-\alpha_{n}$ , as is shown in Figure  \ref{fig:Cn}. It is a $2n$ dimensional representation and naturally gives an order on weights of the fundamental representation. The ordered weights are denoted by $\lambda^{0},\lambda^{1},...,\lambda^{2n-1}$.

\begin{lemma}\label{lem:Cn inner prd}
\begin{equation}
(\lambda^{s},\lambda^{t})=\left\{
                            \begin{array}{ll}
                              \frac{1}{2}, & \hbox{$s=t$;} \\
                              0, & \hbox{$s+t\neq2n-1,s\neq t$;} \\
                              -\frac{1}{2}, & \hbox{$s+t=2n-1$,}
                            \end{array}
                          \right.
\end{equation}where $s,t=0,1,...,2n-1$.
\end{lemma}
Proof: The proof is straightforward.

These inner products are irrelative to the rank $n$.

\begin{lemma}
\begin{equation}
(\lambda^{s},\lambda^{t})=(\lambda^{2n-1-s},\lambda^{2n-1-t}),
\end{equation}where $s,t=0,1,...,2n-1$.
\end{lemma}
Proof: $2n-1-s+2n-1-t=2n-1,2n-1-s=2n-1-t$ if and only if $s+t=2n-1,s=t$. $2n-1-s+2n-1-t\neq2n-1,2n-1-s=2n-1-t$ if and only if $s+t\neq2n-1,s=t$. From the Lemma \ref{lem:Cn inner prd}, the proof is straightforward.

When $\mathbbm{q}=2n-1$ without symmetry breaking, the simple roots $\alpha_{i_{j}}$ in (\ref{braiding formula in NSB}) are $\alpha_{i_{j}}=\left\{
                                                                            \begin{array}{ll}
                                                                              \alpha_{j}, & \hbox{$1\leq j\leq n$;} \\
                                                                              \alpha_{2n-j}, & \hbox{$n<j\leq 2n-1$.}
                                                                            \end{array}
                                                                          \right.$

So\begin{equation}
\mathcal{B}\mathcal{J}_{2n-1}=(-1)^{2n-1}q^{-\frac{1}{2}[(\lambda_{1},\lambda_{2})+\sum_{j<s}(\alpha_{i_{j}},\alpha_{i_{s}})-\sum_{j,\alpha}(\alpha_{i_{j}},\lambda_{\alpha})]}\mathcal{J}_{2n-1}=-q^{\frac{n}{2}+\frac{1}{4}}\mathcal{J}_{2n-1}.
\end{equation}
The minus comes from the fact that the braiding changes each dimension of the thimble $\mathcal{J}_{2n-1}$ into the opposite direction and $\mathcal{J}_{2n-1}$ is odd dimensional.

As the same reason as in the case of $B_{n}$ Lie algebra, it implies that
\begin{equation}\label{equ: fusion matrix as eigenvector in C}
\sum_{c,d}\mathcal{B}_{ab}^{cd}\mathcal{M}_{cd}=-q^{\frac{n}{2}+\frac{1}{4}}\mathcal{M}_{ab}.
\end{equation}

By the same method as in the case of $B_{n}$ Lie algebra, we derive the braiding and fusion matrices for the fundamental representation of $C_{n}$ Lie algebra and omit the proof.

\begin{theorem}
For $\mathbbm{q}\neq 2n-1$,
\begin{equation}
\mathcal{B}\mathcal{J}_{s,\mathbbm{q}-s}=\left\{
                                             \begin{array}{ll}
                                               q^{- \frac{1}{4}}\mathcal{J}_{s,\mathbbm{q}-s}, & \hbox{$s=\mathbbm{q}-s$;} \\
                                               \mathcal{J}_{\mathbbm{q}-s,s}, & \hbox{$s > \mathbbm{q}-s$;} \\
                                               \mathcal{J}_{\mathbbm{q}-s,s}+(q^{-\frac{1}{4}}-q^{\frac{1}{4}})\mathcal{J}_{s,\mathbbm{q}-s}, & \hbox{$s< \mathbbm{q}-s$.}
                                             \end{array}
                                           \right.
\end{equation}
\end{theorem}

\begin{lemma}\label{}
For $\mathbbm{q}=2n-1$,
\begin{equation}
\mathcal{B}\mathcal{J}_{f,2n-1-f}= q^{ \frac{1}{4}}\mathcal{J}_{2n-1-f,f}+\sum_{i=1}^{2n-1-f}\beta_{f,2n-1-f}^{2n-1-f-i,f+i}\mathcal{J}_{2n-1-f-i,f+i}
\end{equation}
\end{lemma}
where $\beta_{f,2n-1-f}^{2n-1-f-i,f+i}$ are unknown constants.

For the fundamental representation of $C_{n}$ Lie algebra, $m=2n-1$, $\gamma=q^{-\frac{1}{4}}$.  From Lemma \ref{Lemma constraint}, \begin{equation}
\mathcal{M}_{a2n-1-a}\mathcal{M}^{a2n-1-a}=\left\{
                                              \begin{array}{ll}
                                              xq^{\frac{n-a}{2}}, & \hbox{$n\leq a\leq 2n-1$;} \\
                                               x^{-1}q^{\frac{n-a-1}{2}}, & \hbox{$0\leq a< n$,}
                                              \end{array}
                                            \right.
\end{equation}where $$x=\mathcal{M}_{nn-1}\mathcal{M}^{nn-1}.$$
 $$C=xq^{-\frac{2n-1}{4}}\hbox{ and  }\beta_{a2n-1-a}^{a2n-1-a}=(q^{-\frac{1}{4}}-q^{\frac{1}{4}})+(q^{\frac{1}{4}}-x^{-2}q^{-\frac{1}{4}})q^{\frac{n-a-1}{2}}.$$

The condition (\ref{norm}) $$|\xi|^{2}=1-\frac{c-c^{-1}}{a-a^{-1}}$$ in Lemma \ref{lem:linear alg} now is equivalent to \begin{equation}\sum_{a,b}(\mathcal{M}_{ab}\mathcal{M}^{ab})=1-\frac{c-c^{-1}}{a-a^{-1}},\end{equation} where $$c=x^{-1}q^{\frac{2n-1}{4}}, a=q^{-\frac{1}{4}}.$$

It follows that $$x=1\hbox{ or }x=-q^{-\frac{1}{2}}.$$

From (\ref{equ: fusion matrix as eigenvector in C}), $x=\mathcal{M}_{nn-1}\mathcal{M}^{nn-1}=-q^{-\frac{1}{2}}$ \begin{equation}
\mathcal{M}_{a2n-1-a}\mathcal{M}^{a2n-1-a}=\left\{
                                              \begin{array}{ll}
                                              -q^{\frac{n-a-1}{2}}, & \hbox{$n\leq a\leq 2n-1$;} \\
                                               -q^{\frac{n-a}{2}}, & \hbox{$0\leq a< n$,}
                                              \end{array}
                                            \right.
\end{equation}
$$C=-q^{-\frac{2n+1}{4}}\hbox{ and  }\beta_{a2n-1-a}^{a2n-1-a}=(q^{-\frac{1}{4}}-q^{\frac{1}{4}})+(q^{\frac{1}{4}}-q^{\frac{3}{4}})q^{\frac{n-a-1}{2}}.$$

From Lemma \ref{prop Laurent ring },
\begin{theorem}
All $2n$ eigenvalues for $\mathcal{B}$ on $V_{2n-1}$ are $$e_{i}=\left\{
                                                                   \begin{array}{ll}
                                                                     q^{-\frac{1}{4}}, & \hbox{$1\leq i\leq n$;} \\
                                                                     -q^{\frac{1}{4}}, & \hbox{$n+1\leq i\leq 2n-1$;} \\
                                                                     -q^{\frac{2n+1}{4}}, & \hbox{$i=2n$.}
                                                                   \end{array}
                                                                 \right.$$
\end{theorem}

Let $y=q^{\frac{1}{4}}$, $a=q^{-\frac{1}{4}}$, $c=-q^{\frac{2n+1}{4}}$, $l=2n$, $m=n$ and $S=\mathcal{B}$ on $V_{2n-1}$, then from Lemma \ref{lem:linear alg},

\begin{equation}\label{equ: skein for C}
\mathcal{B}_{a2n-1-a}^{b2n-1-b}-(\mathcal{B}^{-1})_{a2n-1-a}^{b2n-1-b}=(q^{\frac{1}{4}}-q^{-\frac{1}{4}})(\mathcal{M}_{a 2n-1-a }\mathcal{M}^{b2n-1-b}-\delta_{a}^{b}\delta_{2n-1-a}^{2n-1-b})\hbox{   on }V_{2n-1}.\end{equation}

\begin{theorem}\label{C}
The braiding matrix $\mathcal{B}$ and fusion matrix $M$ corresponding to the fundamental representation of $C_{n}$ Lie algebra satisfy the following conditions:
\begin{enumerate}
  \item $\mathcal{B}_{ac}^{bd}-(\mathcal{B}^{-1})_{ac}^{bd}=(q^{\frac{1}{4}}-q^{-\frac{1}{4}})(\mathcal{M}_{ac} \mathcal{M}^{bd}-\delta_{a}^{b}\delta_{c}^{d})$
  \item $\sum_{a,b}(\mathcal{M}_{ab}\mathcal{M}^{ab})=\frac{q^{-\frac{2n+1}{4}}-q^{\frac{2n+1}{4}}}{q^{\frac{1}{4}}-q^{-\frac{1}{4}}}+1$
  \item $\sum_{c,d}\mathcal{B}^{cd}_{ab}\mathcal{M}_{cd}=-q^{\frac{2n+1}{4}}\mathcal{M}_{ab}$
  \item $\sum_{c,d}(\mathcal{B}^{-1})^{cd}_{ab}\mathcal{M}_{cd}=-q^{-\frac{2n+1}{4}}\mathcal{M}_{ab}.$
\end{enumerate}
\end{theorem}

Under the normalization condition:
\begin{equation}\label{Norm for C ele}
\mathcal{M}_{a2n-1-a}\cdot\mathcal{M}_{2n-1-aa}=1,\quad a=0,1,2,...,2n-1,
\end{equation}
or equivalently,
\begin{equation}\label{Norm for C}
\mathcal{M}_{ab}=\mathcal{M}^{ab},\quad a,b=0,1,2,...,2n-1,
\end{equation} we choose
\begin{equation}\label{equ:fusion C}
\mathcal{M}\left(
                    \begin{array}{c}
                      v_{2n-1} \\
                      v_{2n-2} \\
                      ... \\
                      v_{n} \\
                      v_{n-1} \\
                      ... \\
                      v_{1} \\
                      v_{0} \\
                    \end{array}
                  \right)
=\left(
                     \begin{array}{ccccccccc}
                         &  &  &   &   &   &   & iq^{\frac{-n}{4}} \\
                         &  &  &  &  &  & iq^{\frac{-n+1}{4}} &  \\
                         & &  &  &  & ... &  &  \\
                          &  &  &  &  iq^{-\frac{1}{4}} &  &  &  \\
                         &  &  &  -iq^{\frac{1}{4}} &  &  &  &  \\
                         &  & ... &  &  & 0 &  &  \\
                        & -iq^{\frac{n-1}{4}} &  &  &  &  &  &   \\
                        -iq^{\frac{n}{4}}&  &  &  &  &  &  &    \\
                     \end{array}
                   \right)
\left(
                    \begin{array}{c}
                      v_{2n-1} \\
                      v_{2n-2} \\
                      ... \\
                      v_{n} \\
                      v_{n-1} \\
                      ... \\
                      v_{1} \\
                      v_{0} \\
                    \end{array}
                  \right).
\end{equation}

It should be noticed that Theorem \ref{C} is independent of this normalization.

Now we can use (\ref{equ: skein for C}) to determine the unknown elements $\beta_{a2n-1-a}^{2n-1-b b}$ in the braiding matrix $\mathcal{B}$ on $V_{2n-1}$.

\begin{theorem}For any $a=0,1,2,...2n-2$, $a<b$,
\begin{equation}
\beta_{a2n-1-a}^{2n-1-b b}=\left\{
                              \begin{array}{ll}
                                (q^{\frac{1}{4}}-q^{-\frac{1}{4}})q^{-\frac{1}{4}(a-b)}, & \hbox{$(a-n+\frac{1}{2})(b-n+\frac{1}{2})>0$;} \\
                                (q^{-\frac{1}{4}}-q^{\frac{1}{4}})q^{-\frac{1}{4}(a-b-1)}+\delta_{a}^{2n-1-b}(q^{-\frac{1}{4}}-q^{\frac{1}{4}}), & \hbox{$(a-n+\frac{1}{2})(b-n+\frac{1}{2})<0$.}
                              \end{array}
                            \right.
\end{equation}
\end{theorem}

Thus, the braiding matrix $\mathcal{B}$ is derived and it is independent of the normalization condition (\ref{Norm for C ele}). One can check that when $e^{h}=q^{-\frac{1}{4}}$, this result is the same as that in section 7.3C of \cite{C}, therefore, braiding matrices we derived satisfy Yang-Baxter equation.

\subsection{Fundamental representation of $D_{n}$}

\begin{figure}
\begin{center}
\includegraphics[width=4cm,height=12cm]{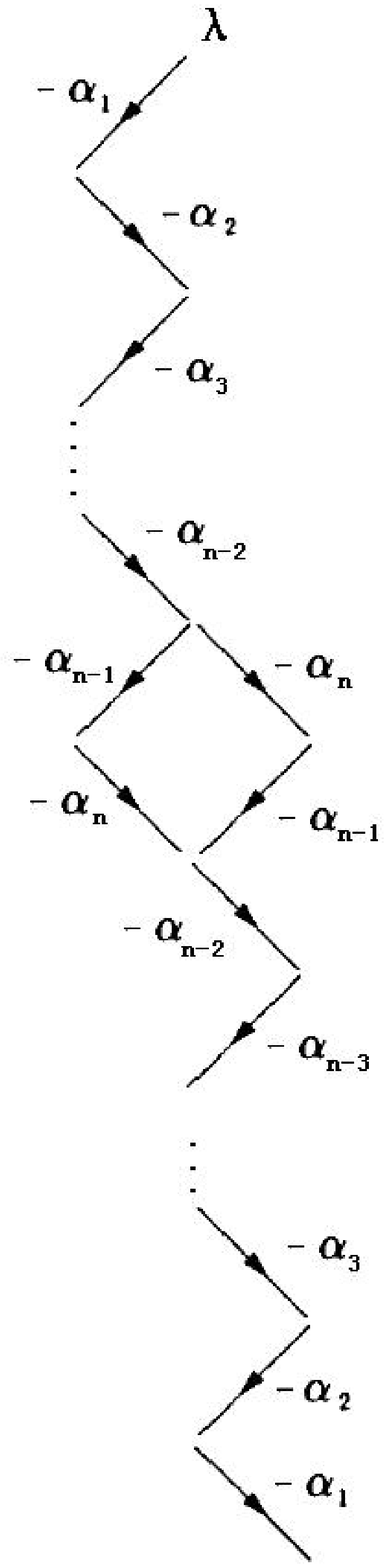}
\caption{\label{fig:Dn}Weights for the fundamental representation of $D_{n}$}
\end{center}
\end{figure}

$D_{n}$ has $n$ simple roots $\alpha_{i}, i=1,2,...,n$. Cartan matrix of $D_{n}$ is
$$\left(
                                                                                        \begin{array}{ccccccc}
                                                                                          2 & -1 & 0 &   &   &   &   \\
                                                                                          -1 & 2 & -1 &   &   &   &   \\
                                                                                          0 & -1 & 2 &   &   &   &   \\
                                                                                            &   &   & ... &   &   &   \\
                                                                                            &   &   &   & 2 & -1 & -1 \\
                                                                                            &   &   &   & -1 & 2 & 0 \\
                                                                                            &   &   &   & -1 & 0 & 2 \\
                                                                                        \end{array}
                                                                                      \right)
.$$
The highest weight of the fundamental representation $V_{\lambda}$ of $D_{n}$ is $\lambda=(1,0,...,0)$. The weights of the fundamental representation are $\lambda,\lambda-\alpha_{1},\lambda-\alpha_{1}-\alpha_{2},...,\lambda-\sum_{i=1}^{n-1}\alpha_{i}, \lambda-\sum_{i=1}^{n-2}\alpha_{i}-\alpha_{n}, \lambda-\sum_{i=1}^{n}\alpha_{i},\lambda-\sum_{i=1}^{n}\alpha_{i}-\alpha_{n-2}, \lambda-\sum_{i=1}^{n}\alpha_{i}-\alpha_{n-2}-\alpha_{n-3},...,\lambda-\sum_{i=1}^{n-2}2\alpha_{i}-\alpha_{n-1}-\alpha_{n}$ , as is shown in Figure  \ref{fig:Dn}. We denote them as  $\lambda^{0},\lambda^{1},...,\lambda^{n-1},\lambda^{{n-1}^{'}},\lambda^{n},...,\lambda^{2n-2}$, where $\lambda^{n-1}=\lambda-\sum_{i=1}^{n-1}\alpha_{i}$ and $\lambda^{{n-1}^{'}}=\lambda-\sum_{i=1}^{n-2}\alpha_{i}-\alpha_{n}$.
For convenience of discussion, we define the order of the weight as follows.
\begin{definition}Order of $a=0,1,...,n-2,n-1,n-1',n,...2n-2$ is defined as
\begin{equation}
o(a)=\left\{
       \begin{array}{ll}
         a, & \hbox{$a=0,1,2,...,n-2,n-1$;} \\
         n, & \hbox{$a=n-1'$;} \\
         a+1, & \hbox{$a=n,n+1,...,2n-2$.}
       \end{array}
     \right.
\end{equation}
\end{definition}

\begin{lemma}
\begin{equation}
(\lambda^{s},\lambda^{t})=\left\{
                            \begin{array}{ll}
                             1, & \hbox{$o(s)+o(t)\neq 2n-1,s=t$;} \\
                              0, & \hbox{$o(s)+o(t)\neq 2n-1,s\neq t$;} \\
                              -1, & \hbox{$o(s)+o(t)=2n-1$,}
                            \end{array}
                          \right.
\end{equation}where $s,t=0,1,2,...,n-2,n-1,n-1^{'},n,n+1,...,2n-2$.
\end{lemma}

Proof: $\lambda^{n-1}+\lambda^{n-1^{'}}= 2\lambda-\sum_{i=1}^{n-2}2\alpha_{i}-\alpha_{n-1}-\alpha_{n}$, $\lambda^{n-1^{'}}+\lambda^{n-1^{'}}\neq 2\lambda-\sum_{i=1}^{n-2}2\alpha_{i}-\alpha_{n-1}-\alpha_{n}$ and $\lambda^{n-1}+\lambda^{n-1}\neq 2\lambda-\sum_{i=1}^{n-2}2\alpha_{i}-\alpha_{n-1}-\alpha_{n}$. If $\lambda^{s}+\lambda^{t}= 2\lambda-\sum_{i=1}^{n-2}2\alpha_{i}-\alpha_{n-1}-\alpha_{n}$, then $s\neq t$. Thus, there are only three cases for the inner products of the weights. The proof is immediate.

Two weights $\lambda^{s}$ and$\lambda^{t}$ are complementary to each other, if $\lambda^{s}+\lambda^{t}=2\lambda-\sum_{i=1}^{n-2}2\alpha_{i}-\alpha_{n-1}-\alpha_{n}$,i.e. $o(s)+o(t)=2n-1$.
These inner products are also irrelative to the rank $n$.
$D_{n}$ also has duality property, the inner products of two weights are the same as the inner products of two complementary weights.
\begin{lemma}
If $o(s)+o(u)=2n-1$ and $o(t)+o(v)=2n-1$, then
\begin{equation}
(\lambda^{s},\lambda^{t})=(\lambda^{u},\lambda^{v}),
\end{equation}where $s,t=0,1,2,...,n-2,n-1,n-1^{'},n,n+1,...,2n-2$.
\end{lemma}
Proof: $o(s)+o(t)=2n-1$ if and only if $o(u)+o(v)=2n-1$. $o(s)+o(t)\neq 2n-1, s\neq t$ if and only if $o(u)+o(v)\neq2n-1, u\neq v$. $o(s)+o(t)\neq 2n-1, s=t$ if and only if $o(u)+o(v)\neq2n-1, u=v$. This completes the proof.

When $\mathbbm{q}=2n-2$ without symmetry breaking, the simple roots $\alpha_{i_{j}}$ in (\ref{braiding formula in NSB}) are  $\alpha_{i_{j}}=\left\{
                                                                            \begin{array}{ll}
                                                                              \alpha_{j}, & \hbox{$1\leq j\leq n$;} \\
                                                                              \alpha_{2n-1-j}, & \hbox{$n<j\leq 2n-2$.}
                                                                            \end{array}
                                                                          \right.$
Thus, \begin{equation}
\mathcal{B}\mathcal{J}_{2n-2}=(-1)^{2n-2}q^{-\frac{1}{2}[(\lambda_{1},\lambda_{2})+\sum_{j<s}(\alpha_{i_{j}},\alpha_{i_{s}})-\sum_{j,\alpha}(\alpha_{i_{j}},\lambda_{\alpha})]}\mathcal{J}_{2n-2}=q^{n-\frac{1}{2}}\mathcal{J}_{2n-2}.
\end{equation}
As the same reason as in the case of $B_{n}$ Lie algebra, it implies that
\begin{equation}\label{equ: fusion matrix as eigenvector in D}
\sum_{c,d}\mathcal{B}_{ab}^{cd}\mathcal{M}_{cd}=q^{n-\frac{1}{2}}\mathcal{M}_{ab}.
\end{equation}
By the same method as in the case of $B_{n}$ Lie algebra, we derive the braiding and fusion matrices for the fundamental representation of $D_{n}$ Lie algebra and omit the proof.

\begin{theorem}
For any $\mathcal{J}_{a,b}\in V_{\mathbbm{q}}$, $\mathbbm{q}\neq 2n-2$,
\begin{equation}
\mathcal{B}\mathcal{J}_{a,b}=\left\{
                                             \begin{array}{ll}
                                               q^{- \frac{1}{2}}\mathcal{J}_{a,b}, & \hbox{$a=b$;} \\
                                               \mathcal{J}_{b,a}, & \hbox{$o(a)>o(b)$;} \\
                                               \mathcal{J}_{b,a}+(q^{-\frac{1}{2}}-q^{\frac{1}{2}})\mathcal{J}_{a,b}, & \hbox{$o(a)<o(b)$.}
                                             \end{array}
                                           \right.
\end{equation}
\end{theorem}

\begin{lemma}\label{}
For any $\mathcal{J}_{a,b}\in V_{\mathbbm{q}}$, $\mathbbm{q}=2n-2$,
\begin{equation}
\mathcal{B}\mathcal{J}_{a,b}=\left\{
                               \begin{array}{ll}
                                 q^{ \frac{1}{2}}\mathcal{J}_{b,a}+\sum_{o(c)<o(b);o(c)+o(d)=2n-1}\beta_{a,b}^{c,d}\mathcal{J}_{c,d}, & \hbox{$a\neq b$;} \\
                                q^{ -\frac{1}{2}}\mathcal{J}_{a,b} , & \hbox{$a=b$.}
                               \end{array}
                             \right.
\end{equation}
\end{lemma}
where $\beta_{a,b}^{c,d}$ are unknown constants.

For the fundamental representation of $D_{n}$ Lie algebra, $m=2n-1$, $\gamma=q^{-\frac{1}{2}}$.  From Lemma \ref{Lemma constraint}, for $a,b$ satisfying $o(a)+o(b)=2n-1$,
\begin{equation}
\mathcal{M}_{ab}\mathcal{M}^{ab}=\left\{
                                              \begin{array}{ll}
                                              xq^{n-o(a)}, & \hbox{$n\leq o(a)\leq 2n-1$;} \\
                                               x^{-1}q^{n-o(a)-1}, & \hbox{$0\leq o(a)< n$,}
                                              \end{array}
                                            \right.
\end{equation}where $$x=\mathcal{M}_{n-1'n-1}\mathcal{M}^{n-1'n-1}.$$
 $$C=xq^{-\frac{2n-1}{2}}\hbox{ and  }\beta_{ab}^{ab}=(q^{-\frac{1}{2}}-q^{\frac{1}{2}})+(q^{\frac{1}{2}}-x^{-2}q^{-\frac{1}{2}})q^{n-o(a)-1}.$$

The condition (\ref{norm}) $$|\xi|^{2}=1-\frac{c-c^{-1}}{a-a^{-1}}$$ in Lemma \ref{lem:linear alg} now is equivalent to \begin{equation}\sum_{a,b}(\mathcal{M}_{ab}\mathcal{M}^{ab})=1-\frac{c-c^{-1}}{a-a^{-1}},\end{equation} where $$c=x^{-1}q^{\frac{2n-1}{2}}, a=q^{-\frac{1}{2}}.$$
It follows that $$x=1\hbox{ or }x=-q^{-1}.$$
From (\ref{equ: fusion matrix as eigenvector in D}),
$x=\mathcal{M}_{n-1'n-1}\mathcal{M}^{n-1'n-1}=1$.

\begin{lemma}\label{Lem tr det for D}
For $\mathbbm{q}=2n-2$, if we choose $\mathcal{J}_{2n-2,0}$, $\mathcal{J}_{2n-3,1}$, ..., $\mathcal{J}_{n,n-2}$, $\mathcal{J}_{n-1,n-1'}$, $\mathcal{J}_{n-1',n-1}$, $\mathcal{J}_{n-2,n}$ ..., $\mathcal{J}_{0,n-2}$, $\mathcal{J}_{n-1,n-1}$, $\mathcal{J}_{n-1,n-1'}$ as a set of basis, $\mathcal{B}$ is a diagonal block matrix $\left(
                                                   \begin{array}{cc}
                                                     \tilde{\mathcal{B}}_{2n\times2n} & 0 \\
                                                     0 & q^{-\frac{1}{2}}I_{2\times2} \\
                                                   \end{array}
                                                 \right)
$ on $V_{2n-2}$,
\begin{equation}
\mathrm{Tr} \tilde{\mathcal{B}}=(n-1)q^{-\frac{1}{2}}-nq^{\frac{1}{2}}+q^{\frac{2n-1}{2}}
\end{equation}
and
\begin{equation}
\mathrm{det} \tilde{\mathcal{B}}=(-1)^{n}q^{n}
\end{equation}
\end{lemma}
If we denote $$\zeta^{t}=(
                                                                          \begin{array}{ccccccc}
                                                                            \zeta^{0} & \zeta^{1}& ...  & \zeta^{n-1} & \zeta^{n-1'} & ... & \zeta^{2n-2} \\
                                                                          \end{array}
                                                                        ),$$
$$ \zeta^{a}=\mathcal{M}^{ab},o(a)+o(b)=2n-1,$$
then \begin{equation}
\mathcal{B}\zeta=C^{-1}\zeta\hbox{ on }V_{2n-1},
\end{equation} i.e. $\zeta$ is an eigenvector of $\tilde{\mathcal{B}}$ on $V_{2n-2}$ with respect to the eigenvalue $C^{-1}=q^{\frac{2n-1}{2}}$.
From Lemma \ref{prop Laurent ring },
\begin{theorem}
all $2n$ eigenvalues for $\tilde{\mathcal{B}}$ on $V_{2n-2}$ are $$e_{i}=\left\{
                                                                   \begin{array}{ll}
                                                                    q^{-\frac{1}{2}}, & \hbox{$1\leq i\leq n-1$;} \\
                                                                     -q^{\frac{1}{2}}, & \hbox{$n\leq i\leq 2n-1$;} \\
                                                                     q^{\frac{2n-1}{2}}, & \hbox{$i=2n$.}
                                                                   \end{array}
                                                                 \right.$$
\end{theorem}

Let $y=q^{\frac{1}{2}}$, $a=q^{-\frac{1}{2}}$, $c=q^{\frac{2n-1}{2}}$, $l=2n$, $m=n-1$ and $S=\tilde{\mathcal{B}}$ on $V_{2n-2}$, then from Lemma \ref{lem:linear alg},
 \begin{equation}\label{equ: skein for D}
\tilde{\mathcal{B}}_{ac}^{bd}-(\tilde{\mathcal{B}}^{-1})_{ac}^{bd}=(q^{\frac{1}{2}}-q^{-\frac{1}{2}})(\mathcal{M}_{a c }\mathcal{M}^{bd}-
\delta_{a}^{b}\delta_{c}^{d})\hbox{   on }V_{2n-2},\end{equation}
$$o(a)+o(c)=o(b)+o(d)=2n-1.$$

\begin{theorem}\label{D}
The braiding matrix $\mathcal{B}$ and fusion matrix $\mathcal{M}$ corresponding to the fundamental representation of $D_{n}$ Lie algebra satisfy the following conditions:
\begin{enumerate}
  \item $\mathcal{B}_{ac}^{bd}-(\mathcal{B}^{-1})_{ac}^{bd}=(q^{\frac{1}{2}}-q^{-\frac{1}{2}})(\mathcal{M}_{ac} \mathcal{M}^{bd}-\delta_{a}^{b}\delta_{c}^{d})$
  \item $\sum_{a,b}(\mathcal{M}_{ab}\mathcal{M}^{ab})=\frac{q^{n-\frac{1}{2}}-q^{-n+\frac{1}{2}}}{q^{\frac{1}{2}}-q^{-\frac{1}{2}}}+1$
  \item $\sum_{c,d}\mathcal{B}^{cd}_{ab}\mathcal{M}_{cd}=q^{n-\frac{1}{2}}\mathcal{M}_{ab}$
  \item $\sum_{c,d}(\mathcal{B}^{-1})^{cd}_{ab}\mathcal{M}_{cd}=q^{-n+\frac{1}{2}}\mathcal{M}_{ab}.$
\end{enumerate}
\end{theorem}

Under the normalization condition:
\begin{equation}\label{Norm for D ele}
\mathcal{M}_{ab}\cdot\mathcal{M}_{ba}=1,\quad a=0,1,2,...,n-1,n-1',...2n-2, o(a)+o(b)=2n-1,
\end{equation}
or equivalently,
\begin{equation}
\mathcal{M}_{ab}=\mathcal{M}^{ab},\quad a,b=0,1,2,...,n-1,n-1',...2n-2,
\end{equation} we choose
\begin{equation}\label{equ:fusion D}
\mathcal{M}\left(
                    \begin{array}{c}
                      v_{2n-2} \\
                      v_{2n-3} \\
                      ... \\
                      v_{n-1'} \\
                      v_{n-1} \\
                      ... \\
                      v_{1} \\
                      v_{0} \\
                    \end{array}
                  \right)
=\left(
                     \begin{array}{ccccccccc}
                         &  &  &   &   &   &   & q^{\frac{-n+1}{2}} \\
                         &  &  &  &  &  & q^{\frac{-n+2}{2}} &  \\
                         & & 0 &  &  & ... &  &  \\
                          &  &  &  & 1 &  &  &  \\
                         &  &  & 1 &  &  &  &  \\
                         &  & ... &  &  & 0 &  &  \\
                        & q^{\frac{n-2}{2}} &  &  &  &  &  &   \\
                        q^{\frac{n-1}{2}}&  &  &  &  &  &  &    \\
                     \end{array}
                   \right)
\left(
                    \begin{array}{c}
                      v_{2n-2} \\
                      v_{2n-3} \\
                      ... \\
                      v_{n-1'} \\
                      v_{n-1} \\
                      ... \\
                      v_{1} \\
                      v_{0} \\
                    \end{array}
                  \right).
\end{equation}

It should be noticed that Theorem \ref{D} is independent of normalization.

Now we can use the formula (\ref{equ: skein for D}) to determine the unknown elements $\beta_{ac}^{d b}$ in the braiding matrix $\tilde{\mathcal{B}}$ on $V_{2n-2}$.

\begin{theorem}For any $a=0,1,2,...2n-3$, $o(a)<o(b)$ and $o(a)+o(c)=o(b)+o(d)=2n-1,$
\begin{equation}
\beta_{ac}^{d b}=\left\{
                              \begin{array}{ll}
                                (q^{\frac{1}{2}}-q^{-\frac{1}{2}})q^{-\frac{1}{2}(o(a)-o(b))}, & \hbox{$(o(a)-n+\frac{1}{2})(o(b)-n+\frac{1}{2})>0$;} \\
                                (q^{\frac{1}{2}}-q^{-\frac{1}{2}})q^{-\frac{1}{2}(o(a)-o(b)+1)}+\delta_{o(a)}^{2n-1-o(b)}(q^{-\frac{1}{2}}-q^{\frac{1}{2}}), & \hbox{$(o(a)-n+\frac{1}{2})(o(b)-n+\frac{1}{2})<0$.}
                              \end{array}
                            \right.
\end{equation}
\end{theorem}

Thus, the braiding matrix $\mathcal{B}$ is derived and it is also independent of the normalization condition (\ref{Norm for D ele}). One can check that when $e^{h}=q^{-\frac{1}{2}}$, this result is the same as that in section 7.3C of \cite{C}, therefore, braiding matrices we derived satisfy Yang-Baxter equation.

\section{Knot invariants for $B_{n},C_{n}$ and $D_{n}$}

Now, we make a brief introduction to Kauffman polynomial. Details can be found in \cite{K}.
\begin{definition}
The equivalence relation on link diagrams generated by Reidemeister move II and III is called regular isotopy; The equivalence relation generated by Reidemeister move I, II and III is called ambient isotopy.
\end{definition}

\begin{definition}
Let $K$ be any oriented link diagram. Writhe of $K$ (twist number of $K$) is defined by the formula $\omega (K)=\sum_{p\in C(K)}\varepsilon(p)$, where $C(K)$ denotes the set of crossings in the diagram of $K$ and $\varepsilon(p)$ denotes the sign of the crossing $p$.
\end{definition}

\begin{figure}
\begin{center}
\includegraphics[width=12cm,height=3cm]{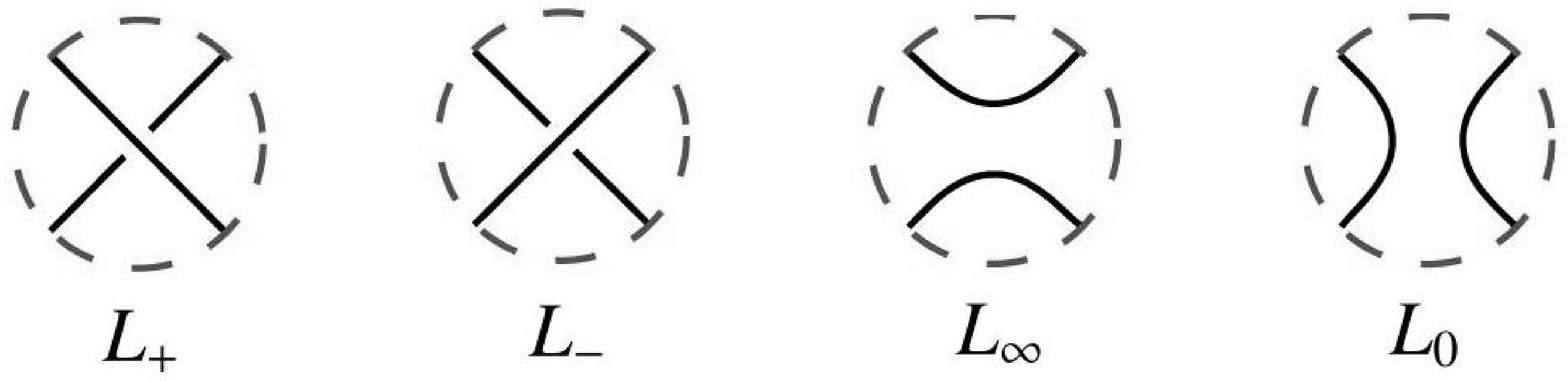}
\caption{\label{fig:4} $L_{+}$, $L_{-}$, $L_{\infty}$ and $L_{0}$}
\end{center}
\end{figure}

\begin{figure}
\begin{center}
\includegraphics[width=8cm,height=3cm]{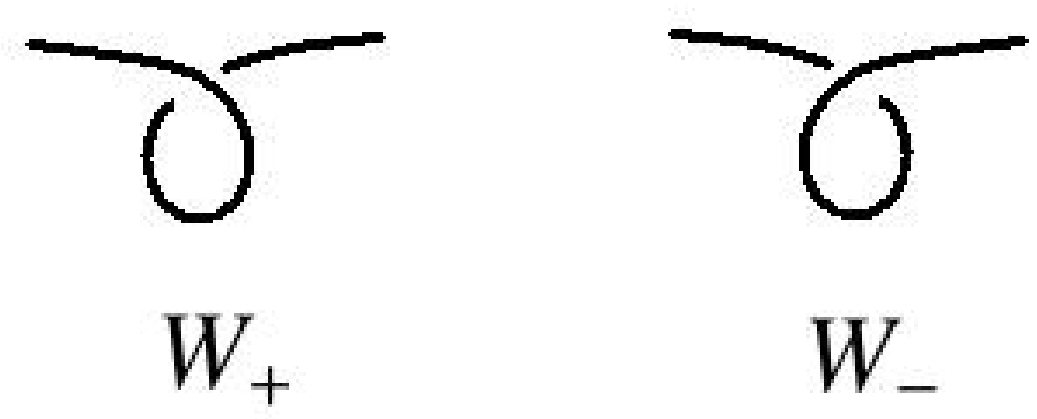}
\caption{\label{fig:2} $W_{+}$ and $W_{-}$}
\end{center}
\end{figure}
Kauffman polynomial is a regular isotopy invariant of unoriented links defined as follows.
\begin{definition}
Kauffman polynomial is a 2-variable polynomial $D_{K}=D_{K}(z,\alpha)$ of unoriented links $K$ satisfying:
\begin{enumerate}
  \item If $K$ and $K'$ are regular isotopy $K\approx K'$, then $D_{K}=D_{K'};$
  \item $D_{L_{+}}-D_{L_{-}}=z(D_{L_{\infty}}-D_{L_{0}});$
  \item $D_{\bigcirc }=((\alpha-\alpha^{-1})/z)+1;$ $D_{W_{+}}=\alpha D_{\mathcal{M}^{-1}};$ $D_{W_{-}}=\alpha^{-1}D_{\mathcal{M}^{-1}},$
\end{enumerate}
where $L_{+}$, $L_{-}$, $L_{\infty}$ and $L_{0}$ are shown in Figure  \ref{fig:4}, and $W_{+}$ and $W_{-}$ are shown in Figure  \ref{fig:2}.
\end{definition}

\begin{figure}
\begin{center}
\includegraphics[width=7.2cm,height=4cm]{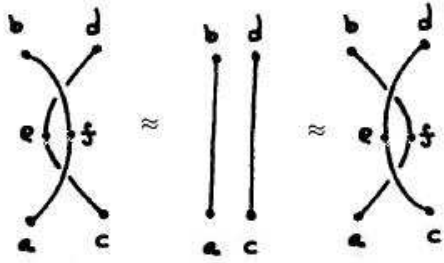}
\caption{\label{fig:VRM} Reidemeister move II}
\end{center}
\end{figure}

Since the operator $\mathcal{B}$ and $\mathcal{B}^{-1}$ are inverse to each other, i.e.
\begin{equation}
\sum_{e,f}(\mathcal{B}^{-1})^{ef}_{ac}\mathcal{B}^{bd}_{ef}=\delta_{a}^{b}\delta_{c}^{d}=\sum_{e,f}\mathcal{B}^{ef}_{ac}(\mathcal{B}^{-1})^{bd}_{ef},
\end{equation}
$<K>$ is clearly invariant under Reidemeister move II (see Figure \ref{fig:VRM}).

Moreover, braiding matrices we derived satisfy Yang-Baxter equation. Thus, $<K>$ is a regular isotopy invariant. $<\bigcirc>$, $<W_{+}>$ and $<W_{-}>$ can be decomposed as \begin{equation}<\bigcirc>=\sum_{a,b}(\mathcal{M}_{ab}\mathcal{M}^{ab}),\end{equation}
\begin{equation}<W_{+}>=\sum_{c,d}\mathcal{B}^{cd}_{ab}\mathcal{M}_{cd}\end{equation}and
\begin{equation}<W_{-}>=\sum_{c,d}(\mathcal{B}^{-1})^{cd}_{ab}\mathcal{M}_{cd},\end{equation}
therefore, Theorem \ref{B}, \ref{C} and \ref{D} show that the knot invariants associated to the simple Lie algebras of type $B_{n}$, $C_{n}$ and $D_{n}$  satisfy the second and third condition in the definition of Kauffman polynomial with $\alpha=q^{n},z=q^{\frac{1}{2}}-q^{-\frac{1}{2}}$; $\alpha=-q^{\frac{2n+1}{4}},z=q^{\frac{1}{4}}-q^{-\frac{1}{4}}$ and $\alpha=q^{\frac{2n-1}{2}},z=q^{\frac{1}{2}}-q^{-\frac{1}{2}}$ respectively.

Thus, for the fundamental representation of simple Lie algebras $B_{n}$, $C_{n}$ and $D_{n}$, the knot invariants $<K>$, defined as vacuum expectations of a quantum mechanics system involving the braiding and fusion operators, are Kauffman polynomial.

\section*{Acknowledgements}

We would like to thank Zhi Chen, A. Losev, E. Witten, Ke Wu and Wenli Yang for very helpful discussions. 
Comments and discussions with Xuexing Lu, Kaiwen Sun, Xiaoyu Jia and Yongjie Wang are also of great help.
This work is partially supported by the National Natural Science Foundation of Grant number
11031005, School of Mathematical Sciences at Capital Normal University and the Wu
Wen Tsun Key Lab of Mathematics of Chinese Academy of Sciences at USTC.


\begin{thebibliography}{90}
{
 \bibitem{J} Jones, Vaughan FR. "A POLYNOMIAL INVARIANT FOR KNOTS VIA VON NEUMANN ALGEBRAS1." (1985).

 \bibitem{K} Kauffman, Louis H. Knots and physics. Vol. 1. World Scientific Publishing Company, 1991.

 \bibitem{W1} Witten, Edward. "Quantum field theory and the Jones polynomial." Communications in Mathematical Physics 121.3 (1989): 351-399.

 \bibitem{GW} Gaiotto, Davide, and Edward Witten. "Knot invariants from four-dimensional gauge theory." arXiv preprint arXiv:1106.4789 (2011).

 \bibitem{FFR} Feigin, Boris, Edward Frenkel, and Nikolai Reshetikhin. "Gaudin model, Bethe ansatz and critical level." Communications in Mathematical Physics 166.1 (1994): 27-62.

 \bibitem{ATY} Awata, Hidetoshi, Akihiro Tsuchiya, and Yasuhiko Yamada. "Integral formulas for the WZNW correlation functions." Nuclear Physics B 365.3 (1991): 680-696.

 \bibitem{W2} Witten, Edward. "Analytic continuation of Chern-Simons theory." Chern-Simons Gauge Theory 20 (2010): 347-446.

 \bibitem{F} Frenkel, Edward. "Free field realizations in representation theory and conformal field theory." Proceedings of the International Congress of Mathematicians. Birkh\"{a}user Basel, 1995.

 \bibitem{C} Chari, Vyjayanthi. "A guide to quantum groups." Cambridge university press, 1995.

 \bibitem{HL} Hu, Sen and Liu, Peng. "HOMFLY polynomial from a generalized Yang-Yang function" arXiv preprint arXiv:1312.1769 (2013).

}
\end{thebibliography}
\end{document}